\documentclass[12pt,reqno]{amsart}
\usepackage{amsmath}
\usepackage{amssymb}
\usepackage{amsthm}
\usepackage{mathrsfs}
\usepackage{latexsym}
\usepackage{amssymb}

\newtheorem{theorem}{Theorem}[section]
\newtheorem{lemma}{Lemma}[section]
\newtheorem{corollary}{Corollary}[section]
\newtheorem{remark}{Remark}[section]
\newtheorem{definition}{Definition}[section]
\newtheorem{proposition}{Proposition}[section]
\newtheorem{example}{Example}[section]
\newtheorem{assumption}{Assumption}[section]
\numberwithin{equation}{section}
\newcommand{\bth}{\begin{theorem}}
\newcommand{\ethe}{\end{theorem}}

\newcommand{\bre}{\begin{remark}}
\newcommand{\ere}{\end{remark}}

\newcommand{\ble}{\begin{lemma}}
\newcommand{\ele}{\end{lemma}}
\newcommand{\bde}{\begin{definition}}
\newcommand{\ede}{\end{definition}}

\newcommand{\bco}{\begin{corollary}}
\newcommand{\eco}{\end{corollary}}

\newcommand{\bpr}{\begin{proposition}}
\newcommand{\epr}{\end{proposition}}

\newcommand{\bexer}{\begin{exercise}}
\newcommand{\eexer}{\end{exercise}}
\newcommand{\breh}{\begin{hint}}
\newcommand{\ereh}{\end{hint}}

\newcommand{\halmos}{\hfill \qed}

\newcommand{\bexam}{\begin{example}}
\newcommand{\eexam}{\end{example}}

\newcommand{\pr} {{\bf Proof.}}

\newcommand{\bfi}{\begin{fig}}
\newcommand{\efi}{\end{fig}}

\newcommand{\beao}{\begin{eqnarray*}}
\newcommand{\eeao}{\end{eqnarray*}\noindent}

\newcommand{\beam}{\begin{eqnarray}}
\newcommand{\eeam}{\end{eqnarray}\noindent}
\newcommand{\E}{\mathbf{E}}
\newcommand{\PP}{\mathbf{P}}

\newcommand{\xto}{x\to\infty}

\newcommand{\vto}{v\to\infty}

\newcommand{\bH}{\overline{H}}
\newcommand{\bF}{\overline{F}}

\newcommand{\bV}{\overline{V}}

\newcommand{\bbr}{{\mathbb R}}

\newcommand{\bbb}{{\mathbb B}}

\newcommand{\bbn}{{\mathbb N}}

\newcommand{\vep}{\varepsilon}

\begin{document}
\title{Random vectors in the presence of a single big jump}

\author[ D.G. Konstantinides, C. D. Passalidis ]{ Dimitrios G. Konstantinides, Charalampos  D. Passalidis}

\address{Dept. of Statistics and Actuarial-Financial Mathematics,
University of the Aegean,
Karlovassi, GR-83 200 Samos, Greece}
\email{konstant@aegean.gr,\;sasd24009@sas.aegean.gr.}

\date{{\small \today}}

\begin{abstract}
Multidimensional distributions with heavy tails have recently attracted the attention of several papers on Applied Probability. About multivariate subexponentiality we can find several approximations, but none of them have been widely established. Having in mind the single big jump principle, and further the multivariate subexponentiality suggested by Samorodnitsky and Sun (2016), we introduce the multivariate long, dominatedly and consistently varying distribution classes. We examine the closure properties of these classes with respect to the product convolution, the scale mixture and the convolution of multivariate distributions. Additionally, for two distributions from the class of multivariate long tailed distributions, we provide necessary and sufficient conditions in order their convolution to belong to the class of multivariate subexponential distributions. Furthermore, we study the multivariate, linear, single big jump principle in finite and in random sums of random vectors, permitting some dependence structures, which contain the independence as a special case. Finally, we present an application on the asymptotic behavior of the entrance probability of the discounted aggregate claims into some rare sets, in a risk model, with a common Poisson counting process, financial factors and independent, identically distributed claims, with common multivariate subexponential distribution.
\end{abstract}

\maketitle
\vspace{3mm}


\section{Introduction} \label{sec.KLP.1}

\subsection{Motivation.}

On the one hand side, it is known that the distributions with heavy tails play a crucial role in applied probability; see for example \cite{embrechts:klueppelberg:mikosch:1997}, \cite{konstantinides:2018}, \cite{asmussen:steffensen:2020} among others. On the other hand side, the dependence modeling does not play a secondary role in the practical applications. Thus, in several branches of applied probability, as in risk theory and in risk management, emerges intensive study of multidimensional models through random vectors; see for example \cite{liu:woo:2014}, \cite{li:2022a}, \cite{liu:shushi:2024}, among others. In these papers we find some multidimensional distribution classes with heavy tails, which represent the multivariate regular variation, symbolically $MRV$, and some of these papers study both $MRV$ and some type of multivariate Gumbel distributions.
 
However, although the multivariate regular variation class is well established; see  \cite{haan:resnick:1981} and \cite{resnick:2007} for introduction and significant treatments, it is somehow restrictive with respect to its distribution range. Indeed, an $MRV$ vector does not contain some moderately heavy-tailed marginal distributions, because of the regular variation of the marginals. Furthermore, in the case of multivariate Gumbel distributions, the restriction comes from the common normalization function, which divides each component and also excludes the dominatedly varying marginal distributions. Similar problems in the $MRV$ are handled through non-standard $MRV$. 

Therefore, these multivariate classes fit better to extreme value theory, but are not sufficient alone for applications on insurance, finance and queuing theory. Thus, recently, several papers have appeared focused on multidimensional distribution classes with heavy tails, with the ultimate aim of the multivariate subexponentiality.

We find at least six different definitions of the multivariate subexponentiality, where we can distinguish two major directions:
\begin{enumerate}
\item
multivariate, non-linear, single, big jump principle: In \cite{konstantinides:passalidis:2024c} were introduced some classes of distributions with heavy tails, in two-dimensional set up, imposing that their subexponentiality satisfies some kind of max-sum equivalence in the two dimensions. The way of definition provides some kind of relation to the dimension of random vectors, as also to the flexibility with respect to the convergence speed of the marginal tails. However, the properties of such distributions are not connected immediately with the corresponding one-dimensional properties. A similar approach can be found in \cite{wang:su:yang:2024}, where some smaller classes of this kind were introduced.

\item
multivariate, linear, single, big jump principle: The papers \cite{cline:resnick:1992}, \cite{omey:2006} and \cite{samorodnitsky:sun:2016} are focused on a kind of single big jump through five different definitions of subexponentiality. These ways of definition contain some kind of 'insensitivity with respect to dimension' in relation to the single big jump principle, while the same time it seems to keep the corresponding one-dimensional properties, in some cases, that is not the case for direction (1). 
\end{enumerate}

In the current paper we follow the direction (2). In first look of these three papers, we find that all of these approaches have a variety of applications. With the first one, \cite{kaleta:ponikowski:2022} was considered the compound Poison process, where the subexponential and convolution equivalent densities are suggested. The second one, examines mostly the asymptotic behavior of randomly stopped sums of random vectors and consequently of also L\'{e}vy processes; see for example \cite{daley:omey:vesilo:2007}, \cite{kapnova:palmowski:2022}. However, these two approaches present a handicap against the third one, where found some good properties, like the preservation of the closure properties of the classes as we step from the one-dimensional to multidimensional set up, or the property for any linear combination of the components of the multivariate subexponential random vector, to have one-dimensional subexponential distribution. In spite of these good properties  from \cite{samorodnitsky:sun:2016}, it did not attract the attention it deserved. Hence, having in mind some applications on risk theory, we suggest some other distribution classes with heavy tails, through a similar way, in order to study the closure properties,  the behavior of the sums of random vectors (with distributions on these classes), and some applications on risk theory.

The paper is organized as follows. In subsection $1.2$ we give some basic concepts of the one-dimensional distributions with heavy tails. In section $2$ after reminding the known multidimensional classes of heavy-tailed distributions, and we introduce three new classes in a multidimensional context, together with some preliminary results. In section 3, we study these new classes, and also the multivariate subexponential class, with respect to some closure properties, as the convolution, the product convolution, the scale mixture and the finite mixture. In section 4, we consider the asymptotic behavior of the sums of random vectors, with the summands belonging to these classes, and with random vectors that are not necessarily independent. We also consider the asymptotic behavior of randomly stopped sums of random vectors. These results highlights the multivariate linear single big jump principle on, the finite, and on the randomly stopped, sums, respectively while such results were available only in $MRV$ case. In section 5, we take into account the asymptotic behavior of the entrance probability of the discounted aggregate claims, into some rare sets, in a multidimensional risk model with a common stochastic discount factor and a Poisson counting process, when the claim distributions belong to the class of multivariate subexponential distributions.

\subsection{Heavy-tailed distributions}

In this section we present some basic definitions of distribution classes with heavy tails, after a short survey of the notations, used further in the whole paper.

For any two positive functions $f,\,g$, we denote by $f(x) \sim g(x)$, $f(x)=o[g(x)]$ and $f(x)=O[g(x)]$, as $\xto$, if 
\beao
\lim_{\xto} \dfrac{f(x)}{g(x)} =1\,, \quad \lim_{\xto} \dfrac{f(x)}{g(x)} =0\,, \quad \limsup_{\xto} \dfrac{f(x)}{g(x)} < \infty\,,
\eeao 
hold respectively, while we write $f(x) \asymp g(x)$, as $\xto$, if both $f(x)=O[g(x)]$ and $g(x)=O[f(x)]$, as $\xto$ hold. If the $d$-dimensional, positive functions ${\bf f},\,{\bf g}$ are of the form ${\bf f}(x\,A),\,{\bf g}(x\,A)$, for some $A \in \bbr^d$, with ${\bf 0} \notin \overline{A}$, where $\overline{A}$ represents the closed case of $A$, then the previous asymptotic relations in one dimension can be generalized in $d$ dimensions, for example ${\bf f}(x\,A)\sim {\bf g}(x\,A)$, as $\xto$, implies $\lim_{\xto} {\bf f}(x\,A)/{\bf g}(x\,A) =1$. For two $d$-dimensional, with $d \in \bbn$, vectors ${\bf X},\,{\bf Y}$, with ${\bf X}=(X_1,\,\ldots,\,X_d)^{\top}$, ${\bf Y}=(Y_1,\,\ldots,\,Y_d)^{\top}$ (where ${\bf X}^{\top}$ denotes the transpose of ${\bf X}$) the operations hold by components, hence they take the form ${\bf X}+{\bf Y}=(X_1+Y_1,\,\ldots,\,X_d+Y_d)^{\top}$ and ${\bf X}\,{\bf Y}=(X_1\,Y_1,\,\ldots,\,X_d\,Y_d)^{\top}$. Further the scalar product is defined for some $\lambda>0$ as  $\lambda\,{\bf X}=(\lambda\,X_1,\,\ldots,\,\lambda\,X_d)^{\top}$, with ${\bf 0}$ the zero vector. For two real numbers $x$, $y$, we denote $x\vee y:=\max\{x,\,y\}$,  $x\wedge y:=\min\{x,\,y\}$, and $\left\lfloor x \right\rfloor$ is the integer part of $x$. 

For any one-dimensional distribution $V$, we denote by $\bV(x) =1- V(x)$ its tail. In the case that the random variable (or, vector) $X$ follows distribution $V$, we denote $X \stackrel{d}{\sim} V$. Let $X_1$ and $X_2$ be two random variables with distributions $V_1$ and $V_2$, and denote the tail of their sum by $\bV_{X_1+X_2}(x) =\PP[X_1 + X_2 >x]$. When $X_1$ and $X_2$ are independent, this tail can be expressed as convolution tail $\overline{V_1*V_2}(x)$. Further, we denote by $V_{X_1 \vee X_2}(x)$ the distribution of the maximum, namely $V_{X_1 \vee X_2}(x)= \PP[X_1\vee X_2\leq x]$.  We say that two distributions $V_1$, $V_2$ are strongly tail equivalent, if holds $\bV_1(x) \sim c\,\bV_2(x)$, as $\xto$, for some $c \in (0,\,\infty)$, while we say that they are weakly tail equivalent, if $\bV_1(x) \asymp \bV_2(x)$, as $\xto$.

For reasons of consistency, we give the following definitions, only in the case where $V$ have support on $\bbr_+=[0,\,\infty)$, although most of these classes can be extended immediately with support on $\bbr=(-\infty,\infty)$. For all distributions with heavy tails we keep in mind that the condition $\bV(x)>0$, for any $x\geq 0$, is valid.

We say that the distribution $V$ has heavy tail, symbolically $V \in \mathcal{K}$, if it holds
\beao
\int_0^{\infty} e^{\vep\,x}\,V(dx) = \infty\,,
\eeao
for any $\vep >0$. A big subclass of the class $\mathcal{K}$ is the class $\mathcal{L}$, of long-tailed distributions. We say that the distribution $V$ has long tail, symbolically $V \in \mathcal{L}$, if it holds
\beao
\lim_{\xto} \dfrac{\bV(x-y)}{\bV(x)} = 1\,,
\eeao
for any (or, equivalently, for some) constant $y >0$. It is well-known that if $V \in \mathcal{L}$, then there exists a function $a\;:\;[0,\,\infty) \to (0,\,\infty)$, named insensitivity function of $V$, such that $a(x) \to \infty$, $a(x) = o(x)$, as $\xto$, and 
\beam \label{eq.AKP.2.b} 
\lim_{\xto} \sup_{|y| \leq a(x)}\left|\dfrac{\bV(x - y)}{\bV(x)} - 1\right| = 0\,,
\eeam 
see for example \cite{foss:korshunov:zachary:2013}, \cite{konstantinides:2018}. The most popular class of heavy-tailed distributions is the class of subexponential distributions. We say  that the distribution $V$ is subexponential, symbolically $V \in \mathcal{S}$, if it holds
\beao
\lim_{\xto} \dfrac{\overline{V^{n*}}(x)}{\bV(x)} = n\,,
\eeao
for any (or, equivalently, for some) integer $n \geq 2$, where the $V^{n*}$ represents the $n$-th fold convolution of distribution $V$ with itself.

These three distribution classes were introduced in \cite{chistyakov:1964}, where it was established that $\mathcal{S} \subsetneq \mathcal{L} \subsetneq \mathcal{K}$; see \cite[Th. 2]{chistyakov:1964}. For applications of subexponentiality in actuarial science and finance we refer to 
\cite{embrechts:klueppelberg:mikosch:1997}, \cite{asmussen:albrecher:2010} and \cite{konstantinides:2018}.

The distribution class of dominatedly varying distributions was introduced in \cite{feller:1969}. We say that the distribution $V$ is dominatedly varying, symbolically $V \in \mathcal{D}$, if it holds
\beao
\limsup_{\xto} \dfrac{\overline{V}(b\,x)}{\bV(x)} < \infty\,,
\eeao
for any (or, equivalently, for some) $b \in (0,\,1)$. It is well-known that $\mathcal{D} \subsetneq \mathcal{K}$,  $\mathcal{D} \not\subseteq \mathcal{S}$, $\mathcal{S} \not\subseteq \mathcal{D}$ and $\mathcal{D} \cap \mathcal{L}= \mathcal{D} \cap \mathcal{S} \neq \emptyset$; see for example in \cite{goldie:1978}. Another even smaller than $\mathcal{D} \cap \mathcal{L}$ is class $\mathcal{C}$ of the consistently varying distribution, introduced by \cite{cline:1994} and \cite{cline:samorodnitsky:1994}. We say that the distribution $V$ is consistently varying, symbolically $V \in C$, if 
\beao
\lim_{b \uparrow 1} \limsup_{\xto} \dfrac{\bV(b\,x)}{\bV(x)}=1\,.
\eeao 

One of the smallest, but important, heavy-tailed distribution classes is that of regularly varying distributions, introduced by \cite{karamata:1933}. We say that the distribution $V$ is regularly varying with index $\alpha \in (0,\,\infty)$, symbolically $V \in \mathcal{R}_{-\alpha}$, if it holds
\beao
\lim_{\xto} \dfrac{\overline{V}(t\,x)}{\bV(x)} =t^{-\alpha}\,,
\eeao
for any $t > 0$. For more properties in this class see in \cite{bingham:goldie:teugels:1987}. It is well-known that 
\beam \label{eq.AKP.2.5} 
\bigcup_{\alpha >0}\mathcal{R}_{-\alpha} \subsetneq \mathcal{C} \subsetneq \mathcal{D}\cap \mathcal{L} \subsetneq \mathcal{S}  \subsetneq \mathcal{L} \subsetneq \mathcal{K}\,.
\eeam
Let close this section with the upper Matuszewska index of a distribution $V$, defined as follows
\beao
J_V^+ = -\lim_{\vto} \dfrac {\log \bV_*(v)}{\log v}\,,
\eeao
where $\bV_*(v)=\liminf_{\xto} \dfrac{\bV(v\,x)}{\bV(x)}$. It is well known that $V \in \mathcal{D}$ if and only if $J_V^+ < \infty$, and if $V \in \mathcal{R}_{-\alpha}$ then $J_V^+ =\alpha$. Further, an important inequality related with the upper Matuszewska index, is the Potter inequality. Namely, if $V\in \mathcal{D}$, then for any $p> J_V^+$ there exist two positive constants $C_1$ and $x_0$, such that for any $x\geq y \geq x_0$ it holds
\beam \label{eq.AKP.2.55}
\dfrac{\bV(y)}{\bV(x)} \leq C_1\,\left( \dfrac xy \right)^p\,. 
\eeam
We refer the reader to \cite[Sec. 2.1]{bingham:goldie:teugels:1987}, \cite[Sec. 2.4]{leipus:siaulys:konstantinides:2023}, for further details about $J_V^+$, and related indexes.

\section{Multivariate heavy-tailed distributions} \label{sec.KP.2}

In this section, we remind some definitions of multivariate distributions with heavy tails and precisely the standard $MRV$ and the class $\mathcal{S_R}$ of multivariate subexponential distributions. Next, we introduce three new distribution classes with heavy tails in multidimensional set up, namely the $\mathcal{D_R}$ of multivariate dominatedly varying distributions, the $\mathcal{L_R}$ of multivariate long-tailed distributions and the $\mathcal{C_R}$ of multivariate consistently varying distributions. Additionally, we provide some initial results for these classes, related with their ordering and some closure properties with respect to their 'tail equivalence'. All the random vectors have distribution support on the $\bbr_+^d:=[0,\,\infty)^d$.

Definition of $MRV$ was given in \cite{haan:resnick:1981}. We say that a random vector ${\bf X}$ follows standard $MRV$, if there exist a Radon measure $\mu$, non-degenerate to zero, and a distribution $V \in \mathcal{R}_{-\alpha}$, with $\alpha \in (0,\,\infty)$, such that it holds
\beam \label{eq.AKP.2.6} 
\lim_{\xto} \dfrac 1{\bV(x)}\,\PP\left[ {\bf X} \in x\,\bbb \right] =\mu(\bbb)\,,
\eeam
for any Borel set $\bbb \subseteq [0,\,\infty]^d$, with ${\bf 0} \notin \overline{\bbb}$, which is $\mu$-continuous. The measure $\mu$, has the property of homogeneity, that means it holds
\beao
\mu\left( \lambda^{1/\alpha}\,\bbb \right)=\lambda^{-1}\,\mu(\bbb)\,,
\eeao
for $\lambda> 0$ and any $\mu$-continuous Borel set $\bbb \subseteq [0,\,\infty]^d$, with ${\bf 0} \notin \overline{\bbb}$. If ${\bf X} \stackrel{d}{\sim} F$ belongs to $MRV$ distributions, then we write $F \in MRV(\alpha,\, V,\,\mu)$.

This class of distributions has been used in the modeling of several problems in applied probability and mathematical statistics; see for example \cite{basrak:davis:mikosch:2002a}, \cite{samorodnitsky:2016},  \cite{chen:liu:2024}, \cite{cheng:konstantinides:wang:2024}, \cite{mikosch:wintenberger:2024} for the importance of this class in heavy-tailed Time Series, Risk theory and Risk management. 

Let us consider a family of sets, represented as 
\beao
\mathscr{R}:=\left\{A \subsetneq \bbr^d \;:\; A \;\text{open,\;increasing},\,A^c\;\text{convex},\, {\bf 0} \notin \overline{A} \right\}\,,
\eeao
where by $A^c$ we denote the complementary set of $A$, and a set $A$ is called 'increasing' if for ${\bf x} \in A$ and any ${\bf a} \in [0,\,\infty)^d$, then it holds  ${\bf x} + {\bf a}\in A $. We should remark, that the set $\mathscr{R}$ represents a cone, with respect to positive scalar multiplication; namely, if $A \in \mathscr{R}$, then for any $\lambda> 0$, we have $\lambda\,A \in \mathscr{R}$.

From \cite[Lem. 4.3]{samorodnitsky:sun:2016} we find that for fixed $A \in \mathscr{R}$, there exists some set $I_A \subsetneq \bbr^d$ such that
\beao 
A=\left\{{\bf x} \in \bbr^d\;:\; {\bf p}^{\top}\,{\bf x}=p_1\,x_1 + \cdots + p_d\,x_d  >1\,,\; \exists \; {\bf p} \in I_A \right\}\,,
\eeao
and further for some fixed $A \in \mathscr{R}$ and ${\bf a} \in \bbr^d=(-\infty,\,\infty)^d$ there exists some $u_1>0$, such that it holds
\beam \label{eq.AKP.2.10} 
(u+u_1)\,A \subsetneq u\,A +{\bf a} \subsetneq (u-u_1)\,A\,,
\eeam
for all $u>u_1$.

Furthermore, from \cite[Lem. 4.5]{samorodnitsky:sun:2016}, we have that for some random vector ${\bf X}$ with distribution $F$, we can define the random variable $Y_A:=\sup\{u\;:\;{\bf X} \in u\,A\}$, following the proper distribution $F_A$, whose tail has the following representation
\beam \label{eq.AKP.2.11} 
\bF_A(x)=\PP[{\bf X} \in x\,A] = \PP\left[\sup_{{\bf p} \in I_A} {\bf p}^{\top}\,{\bf X} >x\right]\,, \qquad \forall \;x>0\,.
\eeam
In such a case, we say that the random vector ${\bf X} \stackrel{d}{\sim} F$ has a multivariate subexponential distribution on $A \in \mathscr{R} $, symbolically $F \in  \mathcal{S}_A$, if $F_A \in  \mathcal{S}$.  Next, we define
\beao
\mathcal{S_{R}}:=\bigcap_{A \in \mathscr{R}} \mathcal{S}_A\,.
\eeao
The class $\mathcal{S_R}$ of multivariate subexponential distributions, introduced in \cite{samorodnitsky:sun:2016}.

Let us notice, that the set 
\beam \label{eq.AKP.2.12} 
\mathcal{G}:=\left\{{\bf x} \;:\; \sum_{i=1}^d l_i\,x_{i} >c \right\}\,,
\eeam
for some $c>0$, $l_1,\,\ldots,\,l_d \geq 0$, with $l_1+\cdots +l_d=1$, belongs to $\mathscr{R}$. Hence, if $F \in \mathcal{S}_\mathscr{R}$ through relation \eqref{eq.AKP.2.11}, it follows that all the non-negative, non-degenerate to zero linear combinations of the components of ${\bf X}$ have one-dimensional subexponential distribution.

Following this definition of subexponentiality, we introduce four distribution classes with heavy tails, in multidimensional set up.

\bde \label{def.AKP.2.1}
Let us consider some fixed $A \in \mathscr{R}$ and some random vector ${\bf X}$ following distribution $F$. We say that ${\bf X}$ has
\begin{enumerate}
\item
multivariate consistently varying distribution on  $A$, symbolically $F \in \mathcal{C}_A$, if $F_A \in \mathcal{C}$.
\item
multivariate dominatedly varying distribution on $A$, symbolically $F \in \mathcal{D}_A$, if $F_A \in \mathcal{D}$.
\item
multivariate long tail on $A$, symbolically $F \in \mathcal{L}_A$, if $F_A \in \mathcal{L}$.
\item
multivariate dominatedly varying, long tail on  $A$, symbolically $F \in (\mathcal{D} \cap \mathcal{L})_A$, if $F_A \in \mathcal{D} \cap \mathcal{L}$.
\end{enumerate}
\ede

For these last classes we introduce the notations 
\beao
\mathcal{C_\mathscr{R}}:=\bigcap_{A\in \mathscr{R}} \mathcal{C}_A\,,\qquad \mathcal{D_\mathscr{R}}:=\bigcap_{A\in \mathscr{R}} \mathcal{D}_A\,,\qquad \mathcal{L_\mathscr{R}}:=\bigcap_{A\in \mathscr{R}} \mathcal{L}_A\,, \quad \mathcal{(D\cap L)_\mathscr{R}}:=\bigcap_{A\in \mathscr{R}} (\mathcal{D}\cap \mathcal{L})_A\,,
\eeao
as the multivariate $\mathcal{C},\,\mathcal{D},\,\mathcal{L},\,\mathcal{D}\cap\mathcal{L}$ over the whole $\mathscr{R}$, respectively.

\bre  \label{rem.AKP.2.1}
It follows directly from the comments in Subsection $1.2$ and by the definitions of these concrete multivariate classes, that, for any fixed $A \in \mathscr{R}$, $(\mathcal{D}\cap \mathcal{L})_A=(\mathcal{D}\cap \mathcal{S})_A$ (in the sense $F \in (\mathcal{D}\cap \mathcal{S})_A$ if $F_A\in (\mathcal{D}\cap \mathcal{S})$) and further we obtain
\beam \label{eq.AKP.2.13} 
\mathcal{C}_A \subsetneq (\mathcal{D}\cap \mathcal{L})_A \subsetneq \mathcal{S}_A \subsetneq \mathcal{L}_A\,.
\eeam
The same inclusions hold for the corresponding classes over $\mathscr{R}$. Now, from relations \eqref{eq.AKP.2.11} and \eqref{eq.AKP.2.12} we find that if $F \in \mathcal{B}_\mathscr{R}$, with $ \mathcal{B} \in \{\mathcal{C},  \mathcal{D},  \mathcal{L}, \mathcal{D} \cap \mathcal{L} \}$, then the non-negative and non-degenerate to zero linear combinations of the components of ${\bf X}$, let say the 
\beao
\sum_{i=1}^d l_i\,X_i\,,
\eeao 
have distribution which belongs to one-dimensional $\mathcal{B}$. This property is significant enough for these distribution classes, since it is related directly to stability and infinite divisibility properties; see for example in \cite[Ch. 2]{samorodnitsky:taqqu:1994}, \cite[Ch.3]{meerschaert:scheffler:2001} respectively.
\ere

\bre \label{rem.AKP.2.2}
Later, in remark \ref{rem.AKP.5.1}, we present some possible forms of set $A$, related with ruin probabilities. The following form $A=\{{\bf x}\;:\;x_i > b_i, \;\forall \; i=1,\,\ldots,\,d \}$ depicts, in some sense, the ruin set, where the $d$-lines of business fail all together. However, we observe that $A \notin \mathscr{R}$, when $d\geq2$. Indeed, for the membership to $\mathscr{R}$, the complement set $A^c$ should be convex. Let assume that ${\bf x},\,{\bf y} \in A^c$, namely they have the form
\beao
\{{\bf x}\;:\;x_i \leq b_i, \;\exists \; i=1,\,\ldots,\,d \}\,,\qquad \{{\bf y}\;:\;y_i \leq b_i, \;\exists \; i=1,\,\ldots,\,d \}\,,
\eeao  
with some integer $d\geq2$, then their convex combination $\lambda\,x_i + (1- \lambda)y_i$, is not necessarily smaller or equal to $b_i$ for some $i=1,...,d$, and any $\lambda \in [0,\,1]$. Hence, set $A^c$ is not convex. For this ruin probability of all the lines of business, the distribution classes from \cite{konstantinides:passalidis:2024c} and \cite{wang:su:yang:2024}, fit better, as they are defined according the common tail behavior. 
\ere

The next result provides insight with respect to ordering of the multivariate classes.

\bpr \label{pr.AKP.2.1}
$MRV(\alpha,\,V,\,\mu) \subsetneq \mathcal{C_\mathscr{R}}$, for $\alpha \in (0,\,\infty)$.
\epr

\pr~
Let $F \in MRV(\alpha,\,V,\,\mu)$, for some $\alpha \in  (0,\,\infty)$. For arbitrarily chosen $A \in \mathscr{R} $, taking into account \cite[proof of Prop. 4.14]{samorodnitsky:sun:2016}, we obtain $\mu(\partial A)=0$, and $\mu(A) \in (0,\,\infty)$. Hence, from relation \eqref{eq.AKP.2.6} we get
\beao
\lim_{\xto} \dfrac 1{\bV(x)}\,\PP\left[Y_A > x\right] =\mu(A)\,.
\eeao
From the fact that $V \in \mathcal{R}_{-\alpha}$, for some $\alpha \in  (0,\,\infty)$ and by the closure property of regular variation with respect to strong tail equivalence; see \cite[Prop. 3.3(i)]{leipus:siaulys:konstantinides:2023}, it is implied that $F_A \in \mathcal{R}_{-\alpha} \subsetneq \mathcal{C}$, as follows from \eqref{eq.AKP.2.5}. Therefore, $F \in \mathcal{C}_A$, and by the arbitrary choice of $A \in \mathscr{R} $, we conclude that $F \in \mathcal{C_\mathscr{R}}$. 
~\halmos

Proposition \ref{pr.AKP.2.1} together with relation \eqref{eq.AKP.2.13} leads to the inclusions
\beam \label{eq.AKP.2.14} 
\bigcup_{0<\alpha < \infty} MRV(\alpha,\,V,\,\mu)  \subsetneq \mathcal{C}_\mathscr{R} \subsetneq \mathcal{(D\cap L)_\mathscr{R}} \subsetneq \mathcal{S_\mathscr{R}} \subsetneq \mathcal{L_\mathscr{R}}\,,
\eeam
which shows, that the ordering of \eqref{eq.AKP.2.5} remains in tact in the frame of multidimensional set up. Note, that when we refer to $A \in \mathscr{R}$, the subscript $A$ can replace the $\mathscr{R}$; see for example \eqref{eq.AKP.2.13} and \eqref{eq.AKP.2.14}.  

\bre \label{rem.AKP.2.3}
The inclusion \eqref{eq.AKP.2.14} is not trivial, in the sense that the classes, greater than $MRV$, do not contain their trivial extensions. We refer the reader to \cite[Sec. 4]{samorodnitsky:sun:2016}, for examples in class $\mathcal{S}_A$, and to \cite[Sec. 4]{konstantinides:liu:passalidis:2025}, for examples in classes $\mathcal{C}_A$ and $(\mathcal{D}\cap \mathcal{L})_A$.
\ere

Next, we assume, that there are two arbitrarily dependent random vectors ${\bf X}$, ${\bf Y}$ with distributions $F$, $G$ respectively, namely $F(x\,A)=\PP[{\bf X} \in x\,A]$, $G(x\,A)=\PP[{\bf Y} \in x\,A]$, for $x \geq 0$, and we explore some closure properties with respect to the asymptotic equivalence of these distributions, in the frame of the newly introduced distribution classes. The following properties comply with the corresponding one-dimensional properties.

\bpr \label{pr.AKP.2.2}
We assume that $A \in \mathscr{R}$, and consider two arbitrarily dependent random vectors ${\bf X}$, ${\bf Y}$ with distributions $F$, $G$ respectively.
\begin{enumerate}
\item[(i)]
If $F \in \mathcal{C}_A$ and
\beam \label{eq.AKP.2.16} 
\lim_{\xto} \dfrac{F (x\,A)}{G (x\,A)}=c\,,
\eeam
for some $c \in (0,\,\infty)$,  then $G \in \mathcal{C}_A$.
\item[(ii)]
If $F \in \mathcal{D}_A$ and
\beam \label{eq.AKP.2.15} 
0 < \liminf_{\xto} \dfrac{F (x\,A)}{G (x\,A)} \leq \limsup_{\xto} \dfrac{F (x\,A)}{G (x\,A)} <\infty\,,
\eeam
then $G\in \mathcal{D}_A$.
\item[(iii)]
If $F \in \mathcal{L}_A$ and
\eqref{eq.AKP.2.16} is valid for some $c \in (0,\,\infty)$, then $G \in \mathcal{L}_A$.
\item[(iv)]
If $F \in (\mathcal{D}\cap \mathcal{L})_A$ and
\eqref{eq.AKP.2.16} is valid for some $c \in (0,\,\infty)$, then $G \in (\mathcal{D}\cap \mathcal{L})_A$.
\end{enumerate}
\epr 

\pr~
Let us consider the random variable $Y_A':=\sup\{u\;:\; {\bf Y} \in u\,A \}$, with distribution $F_A'$.
\begin{enumerate}
\item[(i)]
By relation \eqref{eq.AKP.2.16}, since $\bF_A(x) \sim c\,\bF_A'(x)$, as $\xto$, and by the closure property of class $\mathcal{C}$ with respect to strong tail equivalence; see for example \cite[Prop. 3.5(i)]{leipus:siaulys:konstantinides:2023}, we find that $F_A' \in \mathcal{C}$ and therefore  $G \in \mathcal{C}_A$.
\item[(ii)]
By relation \eqref{eq.AKP.2.15}, with the help of relation \eqref{eq.AKP.2.11}, we obtain $\bF_A(x) \asymp \bF_A'(x)$, as $\xto$, and by closure property of class $\mathcal{D}$ with respect to weak tail equivalence; see for example \cite[Prop. 3.7(i)]{leipus:siaulys:konstantinides:2023}, we find that $F_A' \in \mathcal{D}$ and therefore  $G \in \mathcal{D}_A$.
\item[(iii)]
By relation \eqref{eq.AKP.2.16}, we obtain $\bF_A(x) \sim c\,\bF_A'(x)$, as $\xto$, and by the closure property of class $\mathcal{L}$ with respect to strong tail equivalence; see for example \cite[Prop. 3.9(i)]{leipus:siaulys:konstantinides:2023}, we find that $F_A' \in \mathcal{L}$ and therefore  $G \in \mathcal{L}_A$.
\item[(iv)]
It follows immediately by the two previous parts and the fact that \eqref{eq.AKP.2.16} is included in \eqref{eq.AKP.2.15}.~\halmos
\end{enumerate}

The way we introduced class $\mathcal{L}_A$, allows us to find an insensitivity function for the asymptotic behavior of random vectors belonging to this class. Namely, if $F \in \mathcal{L}_A$, then $F_A \in \mathcal{L}$ and hence there exists some insensitivity function $a\;:\;[0,\,\infty) \to (0,\,\infty)$, such that \eqref{eq.AKP.2.b} holds, and hence by \eqref{eq.AKP.2.11} we obtain
\beao
\lim_{\xto} \sup_{|y| \leq a(x)} \left| \dfrac{\PP\left[{\bf X} \in (x - y)\,A\right]}{\PP\left[{\bf X} \in x\,A\right]} -1 \right| = 0 \,.
\eeao
In next result we find another property of distribution class $\mathcal{L}_A$.

\bpr \label{pr.AKP.2.3}
If we assume that $A \in \mathscr{R}$ is a fixed set and $F \in \mathcal{L}_A$, then it holds
\beam \label{eq.AKP.2.18} 
\lim_{\xto} \dfrac{F (x\,A+{\bf a})}{F (x\,A)}=1\,,
\eeam
for any ${\bf a} \in \bbr^d$.
\epr 

\pr~
Let us consider initially ${\bf a}={\bf 0}$, then relation \eqref{eq.AKP.2.18} follows immediately. Now, we consider that ${\bf a} \in \bbr^d \setminus \{{\bf 0}\} $. Then by relation \eqref{eq.AKP.2.10}, there exists some $u_1>0$, such that
\beam \label{eq.AKP.2.19} 
\lim_{\xto} \dfrac{F(x\,A+{\bf a})}{F (x\,A)}\leq \lim_{\xto} \dfrac{F [(x-u_1)\,A]}{F (x\,A)}=\lim_{\xto} \dfrac{\PP\left[Y_A > x-u_1\right]}{\PP\left[Y_A > x\right]}=1\,,
\eeam
where obviously the inequality $x>u_1$ holds, since $\xto$, and in the last step we used the property of class $\mathcal{L}$. From the other side, using the first inclusion of \eqref{eq.AKP.2.10} and the property of  class $\mathcal{L}$ we obtain
\beam \label{eq.AKP.2.20} 
\lim_{\xto} \dfrac{F (x\,A+{\bf a})}{F (x\,A)}\geq \lim_{\xto} \dfrac{F [(x+u_1)\,A]}{F (x\,A)}=\lim_{\xto} \dfrac{\PP\left[Y_A > x+u_1\right]}{\PP\left[Y_A > x\right]}=1\,,
\eeam
and so through relations \eqref{eq.AKP.2.19} and \eqref{eq.AKP.2.20} we conclude \eqref{eq.AKP.2.18}.  
~\halmos

Note that \eqref{eq.AKP.2.18} does not imply necessarily membership in $\mathcal{L}_A$.
Next, we generalize the \cite[Lem. 4.9]{samorodnitsky:sun:2016}, and the proof is very similar. Let us adopt the following notation
\beam \label{eq.AKP.2.21}
&&{\bf X}^{(i)}=\left(X_1^{(i)},\,\ldots,\, X_d^{(i)}\right)^{\top},\,{\bf S}_n={\bf X}^{(1)}+\cdots+{\bf X}^{(n)},\,Y_A^{(i)}:= \sup\left\{u\;:\; {\bf X}^{(i)} \in u\,A\right\}\,, \\[2mm] \notag
&&\bF_{Y_A^{(1)}+\cdots+Y_A^{(n)}}(x)= \PP\left[Y_A^{(1)}+\cdots+Y_A^{(n)} > x\right],\; F_{{\bf S}_{n}}(x A):= \PP\left[{\bf X}^{(1)}+\cdots+{\bf X}^{(n)} \in x A\right]\,,
\eeam
for $i=1,\,\ldots,\,n$, with $n \in \bbn$, where the random vectors ${\bf X}^{(1)},\,\ldots,\,{\bf X}^{(n)}$ are arbitrarily dependent and not necessarily identically distributed.

\bpr \label{pr.AKP.2.4}
For any fixed set $A \in \mathscr{R}$, it holds
\beam \label{eq.AKP.2.22} 
\bF_{Y_A^{(1)}+\cdots+Y_A^{(n)}}(x) \geq F_{{\bf S}_n}(x\,A)\,,
\eeam
for any $n \in \bbn$, $x \geq 0$.
\epr 

\pr~
From relation \eqref{eq.AKP.2.11} for $x\geq 0$, we obtain
\beao
F_{{\bf S}_n}(x\,A)&=& \PP\left[{\bf X}^{(1)}+\cdots+{\bf X}^{(n)} \in x A\right]\\[2mm]
&=& \PP\left[\sup_{{\bf p} \in I_A} {\bf p}^{\top}\,\left({\bf X}^{(1)}+\cdots+{\bf X}^{(n)}\right)> x \right]\\[2mm]
&\leq & \PP\left[\sup_{{\bf p} \in I_A} {\bf p}^{\top}\,{\bf X}^{(1)}+\cdots+\sup_{{\bf p} \in I_A} {\bf p}^{\top}\,{\bf X}^{(n)}> x \right]\\[2mm]
&=&\PP\left[Y_A^{(1)}+\cdots+Y_A^{(n)} > x\right]\,,
\eeao
and now from \eqref{eq.AKP.2.21} we derive relation \eqref{eq.AKP.2.22}.  
~\halmos

\section{Closure properties} \label{sec.KP.3}

It is well known that the closure properties of heavy-tailed distributions can be useful in applications; see for example \cite{leipus:siaulys:konstantinides:2023} for the one-dimensional distributions and  \cite{hult:lindskog:2006b}, \cite{fougeres:mercadier:2012}, \cite{das:fasenhartmann:2023} for the closure properties of $MRV$. In this section, we provide some closure properties of the multidimensional distribution classes, previously introduced, as also for multivariate subexponential distributions.

\subsection {Product convolution and scale mixture.}

Here, we study the closure properties of class $\mathcal{D}_A$ with respect to product convolution of random vectors, and the closure properties of classes $\mathcal{L}_A$, $(\mathcal{D} \cap \mathcal{L})_A$, $\mathcal{C}_A$, $\mathcal{S}_A$ with respect to scale mixture. 

In case of two random vectors ${\bf \Theta}=(\Theta_{1},\,\ldots,\,\Theta_{d})^{\top}$ and ${\bf X}=(X_1,\,\ldots,\,X_d)^{\top}$, the Hadamard product of random vectors is defined component-wise, as $
{\bf \Theta}\,{\bf X}=(\Theta_{1}\,X_1,\,\ldots,\,\Theta_{d}\,X_d)^{\top}$. To avoid the trivial cases, we assume that the $\Theta_{1},\,\ldots,\,\Theta_{d}$ are non-degenerate to zero. For the previous product we define the random variable
\beao
Z_A:=\sup \{u\;:\; {\bf \Theta}\,{\bf X} \in u\,A\}\,,
\eeao
with distribution $H_A$. In the special case of ${\bf \Theta}= \Theta\,{\bf 1}$, we find the scale mixture $\Theta\,{\bf X}=(\Theta\,X_1,\,\ldots,\,\Theta\,X_d)^{\top}$, for which we define the following random variable
\beam \label{eq.AKP.3.24} 
M_A:=\sup \{u\;:\; \Theta\,{\bf X} \in u\,A\}\,,
\eeam 
with distribution $G_A$. Let us remind that all the random vectors have non-negative components. The next result shows a stability property of class $\mathcal{D}_A$, with respect to product convolution.

\bth \label{th.AKP.3.1}
Let $A \in \mathscr{R}$ be a fixed set, and ${\bf X} \stackrel{d}{\sim} F\in\mathcal{D}_A$. If ${\bf \Theta}$ and ${\bf X}$ are independent, and $\E[\Theta_i^p]< \infty$, for some $p> J_{F_A}^+$, and any $i=1,\,\ldots,\,d$, then for the product convolution ${\bf\Theta\,\bf X} \stackrel{d}{\sim} H$, it holds $H \in \mathcal{D}_A$, and further we obtain $H(x\,A) \asymp F(x\,A)$, as $\xto$.
\ethe

\pr~
We shall show that $\bH_A(x) \asymp \bF_A(x)$. On the one hand, for $\check{\Theta}:=\bigwedge_{i=1}^d \Theta_i$, we find for $x>0$
\beao
\bH_A(x)&=&\PP\left[Z_A >x\right] =\PP\left[\sup_{{\bf p} \in I_A} {\bf p}^{\top}\,{\bf \Theta}\,{\bf X} > x\right] \\[2mm]
&\geq& \PP\left[\sup_{{\bf p} \in I_A} {\bf p}^{\top}\,\check{\Theta}\,{\bf X} > x\right] = \PP\left[ \check{\Theta}\,Y_{A} > x\right]\,.
\eeao

On the other hand, following similar arguments for $\widehat{\Theta}:=\bigvee_{i=1}^d \Theta_i$, for any $x>0$, we obtain
\beao
\bH_A(x)&=&\PP\left[\sup_{{\bf p} \in I_A} {\bf p}^{\top}\,{\bf \Theta}\,{\bf X} > x\right] \\[2mm]
&\leq& \PP\left[\sup_{{\bf p} \in I_A} {\bf p}^{\top}\,\widehat{\Theta}\,{\bf X} > x\right]= \PP\left[ \widehat{\Theta}\,Y_{A} > x\right]\,.
\eeao
Therefore, for any $x>0$ it holds
\beam \label{eq.AKP.3.24.1} 
\PP\left[ \check{\Theta}\,Y_{A} > x\right] \leq \bH_A(x) \leq \PP\left[ \widehat{\Theta}\,Y_{A} > x\right]\,.
\eeam
Further, by \cite[Th. 3.3(iv)]{cline:samorodnitsky:1994}, because of the moment condition on the component of ${\bf \Theta}$ and the independence of ${\bf \Theta}$ and ${\bf X}$,  we obtain as $\xto$ :  
\beam \label{eq.AKP.3.24.2} 
\PP\left[ \check{\Theta}\,Y_{A} > x\right] \asymp \bF_A(x)\,, \qquad  \PP\left[ \widehat{\Theta}\,Y_{A} > x\right] \asymp \bF_A(x)\,.
\eeam
From \eqref{eq.AKP.3.24.1} and \eqref{eq.AKP.3.24.2} we find
\beao
\limsup_{\xto} \dfrac{\bH_A(x)}{\bF_A(x)} \leq \limsup_{\xto} \dfrac{\PP\left[ \widehat{\Theta}\,Y_{A} > x\right]}{\bF_A(x)} < \infty\,,
\eeao
and 
\beao
\liminf_{\xto} \dfrac{\bH_A(x)}{\bF_A(x)} \geq \liminf_{\xto} \dfrac{\PP\left[ \check{\Theta}\,Y_{A} > x\right]}{\bF_A(x)} >0\,.
\eeao
From these two last relations we obtain $\bH_A(x) \asymp \bF_A(x)$, and consequently $H(x\,A) \asymp F(x\,A)$, which, since $F\in \mathcal{D}_A$, from Proposition \ref{pr.AKP.2.2} (ii) follows the inclusion $H \in \mathcal{D}_A$.
~\halmos

Now we restrict ourselves to the case of scale mixture. The following assumption is an extension of a popular assumption, in closure with respect to product in one-dimensional heavy-tailed distributions; see for example \cite{tang:2006}, \cite{konstantinides:leipus:siaulys:2022}, \cite{konstantinides:passalidis:2024b}.

\begin{assumption} \label{ass.AKP.3.1}
It holds
\beam \label{eq.AKP.3.25} 
\PP[\Theta> c\,x]=o\left(\PP\left[M_A>x\right]\right)\,,
\eeam
as $\xto$, for each $c>0$.
\end{assumption}

By relations \eqref{eq.AKP.3.24} and \eqref{eq.AKP.3.25}, we obtain that this is equivalent to the asymptotic relation 
\beao
\PP[\Theta > c\,x]= o\left( \PP\left[\Theta\,{\bf X} \in x\,A\right] \right)= o\left( \PP\left[\sup_{{\bf p} \in I_A} {\bf p}^{\top}\,\Theta\,{\bf X} > x\right] \right)= o\left( \PP\left[\Theta\,Y_A > x\right] \right)\,,
\eeao 
as $\xto$. Furthermore, if the random variable $\Theta$ has a support with a finite right endpoint, and $Y_A$ has a support with infinite right endpoint, then Assumption \ref{ass.AKP.3.1} is true immediately. The scale mixture is omnipresent in several applications, although it seems quite restrictive in relation to product convolution, for example see \cite{li:sun:2009}, \cite{zhu:li:2012}.

Before giving the main theorem for the scale mixtures, we note that the distribution of $\Theta$ is denoted by $Q$, while the set of all positive discontinuity points of $Q$ is denoted by $D(Q)$.

\bth \label{th.AKP.3.1b}
Let $A \in \mathscr{R}$ be a fixed set, and we assume that the distribution of ${\bf X}$ is $F$.
\begin{enumerate}

\item[(i)]
If $\Theta$ and ${\bf X}$ are independent, and they satisfy Assumption \ref{ass.AKP.3.1} with $F \in \mathcal{L}_A$, then for the scale mixture $\Theta\,{\bf X}$ with distribution $G$, it holds $G \in \mathcal{L}_A$.
 
\item[(ii)]
If the assumptions of the previous part hold, with the restriction $F \in (\mathcal{D} \cap \mathcal{L})_A$, then for the scale mixture $\Theta\,{\bf X}$ with distribution $G$, it holds $G \in  (\mathcal{D} \cap \mathcal{L})_A$.
 
\item[(iii)]
If the assumptions of the previous part hold, with the restriction $F \in \mathcal{C}_A$, then for the scale mixture $\Theta\,{\bf X}$ with distribution $G$, it holds $G \in  \mathcal{C}_A$.

\item[(iv)]
Let $\Theta$ be independent of ${\bf X}$, with $F \in \mathcal{S}_A$. Then for the scale mixture $\Theta\,{\bf X}$ with distribution $G$, it holds $G \in \mathcal{S}_A$, if and only if either $D(Q) = \emptyset$ or $D(Q) \neq \emptyset$ and it holds
\beao
\overline{Q}\left(\dfrac xd \right) - \overline{Q}\left(\dfrac {x+1}d \right) = o(\PP\left[M_A >x\right])\,,
\eeao
as $\xto$, for all $d \in D(Q) $.
\end{enumerate}
\ethe

\pr~
Let denote by $\mathcal{B}$ a general distribution class from the set $\{\mathcal{C},\,(\mathcal{D}\cap \mathcal{L}),\,\mathcal{S},\,\mathcal{L}\}$. The inclusion $G \in \mathcal{B}_A $ is equivalent to $M_A \stackrel{d}{\sim} G_A \in \mathcal{B}$. However, from relation \eqref{eq.AKP.2.11} and \eqref{eq.AKP.3.24}  we obtain $\PP[M_A >x]= \PP[\Theta\,Y_A >x]$, for $x>0$. Hence, it is enough to show that the distribution of $\Theta\,Y_A$ belongs to class $\mathcal{B}$. Since $F \in \mathcal{B}_A$, thus we obtain $Y_A \stackrel{d}{\sim} F_A \in \mathcal{B}$ and from the fact that the $\Theta$, $Y_A$ are independent, non-negative random variables, because of the assumptions in points $(i) - (iv)$, the result follows directly from the following collective statement

\begin{enumerate}

\item[(i)]
For the class $\mathcal{L}$, the result is implies by \cite[Th. 2.2(iii)]{cline:samorodnitsky:1994}. 

\item[(ii)]
For the class $(\mathcal{D}\cap \mathcal{L})$, the result is implied by \cite[Cor. 5.2(b)]{leipus:siaulys:konstantinides:2023}.

\item[(iii)]
For the class $\mathcal{C}$, the result is implied by \cite[Th. 3.4(ii)]{cline:samorodnitsky:1994}.

\item[(iv)]
For the class $\mathcal{S}$, the result is implied by \cite[Th. 1.3]{xu:cheng:wang:cheng:2018}.
~\halmos
\end{enumerate}

\subsection{Convolution and finite mixture.}

In this subsection we study the closure of some distribution classes with respect to convolution and finite mixture. The motivation to study these two properties together stems from \cite{embrechts:goldie:1980}, \cite{leipus:siaulys:2020} where the closure of one-dimensional subexponential distributions with respect to convolution was examined. In the case of independence, we simplify the notation from \eqref{eq.AKP.2.21} as follows: let  ${\bf X}^{(1)}:=\left(X_1^{(1)},\,\ldots,\,X_d^{(1)}\right)^{\top} $ and  ${\bf X}^{(2)}:=\left(X_1^{(2)},\,\ldots,\,X_d^{(2)}\right)^{\top}$ be independent random vectors, with distributions $F_1$  and $F_2$ respectively, then we define their convolutions as follows $F_1*F_2(x\,A) =\PP\left[{\bf X}^{(1)} + {\bf X}^{(2)} \in x\,A\right]$, and we put $Y_A^{(1)}$ and $Y_A^{(2)}$, as in comments before Proposition \ref{pr.AKP.2.4}. We have to note that for a multivariate distribution class $\mathcal{B}_A$, with $A \in \mathscr{R}$, we say that the convolution $F_1*F_2 \in \mathcal{B}_A$, if for the random variable 
\beam \label{eq.AKP.3.32} 
Y_A^* :=\sup\left\{ u\;:\;{\bf X}^{(1)} + {\bf X}^{(2)} \in u\,A\right\}\,,
\eeam
with distribution $F_A^*$, it holds that $F_A^* \in \mathcal{B}$. We denote $\overline{F_A^{(1)}*F_A^{(2)}}(x) =  \PP\left[Y_A^{(1)} + Y_A^{(2)} >x\right]$. 

For the closure with respect to finite mixture it is enough to show
\beam \label{eq.AKP.3.33} 
p\,F_1(x\,A) +(1-p)\,F_2(x\,A)  \in \mathcal{B}_A\,,
\eeam
for any $p \in (0,\,1)$, which is equivalent to
\beam \label{eq.AKP.3.34} 
p\,\overline{F_A^{(1)}}(x) +(1-p)\,\overline{F_A^{(2)}}(x)  \in \mathcal{B}\,,
\eeam
where by $F_A^{(i)}$ is  denoted the distribution tail of the $Y_A^{(i)}$ for $i=1,\,2$.

\bpr \label{pr.AKP.3.1}
We assume that $A \in \mathscr{R}$ is a fixed set, $F_1 \in \mathcal{B}_A$, and

a) $F_2 \in \mathcal{B}_A$, then \eqref{eq.AKP.3.33} holds, for any $p \in (0,\,1)$, where $\mathcal{B}_A \in \left\{\mathcal{C}_A, \mathcal{D}_A, \mathcal{L}_A, (\mathcal{D} \cap \mathcal{L})_A \right\}$.

b) $F_2(x\,A) = o\left[F_1(x\,A)\right]$, as $\xto$. 
Then relation \eqref{eq.AKP.3.33} holds, for any $p \in (0,\,1)$, where $\mathcal{B}_A \in \{\mathcal{C}_A,\,\mathcal{L}_A\}$.
\epr

\pr~
In both cases, it follows by the equivalence of \eqref{eq.AKP.3.33} and  \eqref{eq.AKP.3.34}, in combination with \cite[Prop. 3.5(iii)]{leipus:siaulys:konstantinides:2023}  for $\mathcal{C}$, with \cite[Prop. 3.7(iv)]{leipus:siaulys:konstantinides:2023} for $\mathcal{D}$,  with \cite[Prop. 3.9(iii)]{leipus:siaulys:konstantinides:2023} for $\mathcal{L}$ and with \cite[Prop. 3.10(iii)]{leipus:siaulys:konstantinides:2023} for $\mathcal{D}\cap \mathcal{L}$.
~\halmos

Next statement establishes the closure property of distribution classes $\mathcal{D}_A$,  $(\mathcal{D} \cap \mathcal{L})_A $ and $\mathcal{C}_A$, together with some kind of weak or strong max-sum equivalence, in multidimensional set up.

\bth \label{th.AKP.3.2}
We assume that $A \in \mathscr{R}$ is a fixed set and ${\bf X}^{(1)}$, ${\bf X}^{(2)}$ are independent random vectors with distributions $F_1$, $F_2$ respectively.
\begin{enumerate}

\item[(i)]
If $F_1,\,F_2 \in \mathcal{D}_A$, then 
\beam \label{eq.AKP.3.35} 
F_1*F_2(x\,A) \asymp F_1(x\,A) +F_2(x\,A)\,,
\eeam
as $\xto$, and further it holds $F_1*F_2 \in \mathcal{D}_A$. 

\item[(ii)]
If $F_1,\,F_2 \in \mathcal{(\mathcal{D} \cap \mathcal{L})}_A$, then 
\beam \label{eq.AKP.3.36} 
F_1*F_2(x\,A) \sim  F_1(x\,A) +F_2(x\,A)\,,
\eeam
as $\xto$, and further it holds $F_1*F_2 \in \mathcal{(\mathcal{D} \cap \mathcal{L})}_A$. 

\item[(iii)]
If $F_1,\,F_2 \in \mathcal{C}_A$, then  it holds \eqref{eq.AKP.3.36} and further it holds $F_1*F_2 \in \mathcal{C}_A$. 
\end{enumerate} 
\ethe

\pr~
\begin{enumerate}

\item[(i)]
By \cite[Prop. 2.1]{cai:tang:2004}, see also in \cite[Prop. 3.7 (ii)]{leipus:siaulys:konstantinides:2023}, we obtain 
\beam \label{eq.AKP.3.37} 
\overline{F_A^{(1)}*F_A^{(2)}}(x) \asymp \overline{F_A^{(1)}}(x)+\overline{F_A^{(2)}}(x)\,,
\eeam
as $\xto$,  since $F_A^{(i)}\in \mathcal{D}$ for $i=1,\,2$. Further this provides that $F_A^{(1)}*F_A^{(2)} \in \mathcal{D}$, which is implied
by relation \eqref{eq.AKP.3.37}, from the closure of class $\mathcal{D}$ with respect to weak tail equivalence and by the inclusion $\overline{F_A^{(1)}}+\overline{F_A^{(2)}} \in \mathcal{D}$. Indeed, this last comes from the definition of class $\mathcal{D}$, since it holds
\beao
\limsup_{\xto} \dfrac{\overline{F_A^{(1)}}(b\,x)+\overline{F_A^{(2)}}(b\,x)}{\overline{F_A^{(1)}}(x)+\overline{F_A^{(2)}}(x)} \leq \max\left\{\limsup_{\xto} \dfrac{\overline{F_A^{(1)}}(b\,x)}{\overline{F_A^{(1)}}(x)}\,,\;\limsup_{\xto} \dfrac{\overline{F_A^{(2)}}(b\,x)}{\overline{F_A^{(2)}}(x)} \right\} < \infty\,,
\eeao
for any $b \in (0,\,1)$. 

Now, by relation \eqref{eq.AKP.3.37} and Proposition \ref{pr.AKP.2.4} we find the upper bound
\beam \label{eq.AKP.3.38} 
\limsup_{\xto} \dfrac{\PP\left[{\bf X}^{(1)}+{\bf X}^{(2)} \in x\,A\right]}{\PP\left[{\bf X}^{(1)} \in x\,A\right]+\PP\left[{\bf X}^{(2)} \in x\,A\right]} \leq \limsup_{\xto} \dfrac{\overline{F_A^{(1)}*F_A^{(2)}}(x)}{\overline{F_A^{(1)}}(x)+\overline{F_A^{(2)}}(x)} <\infty\,.
\eeam

For the lower bound, because we have non-negative random vectors and the set $A$ is by definition increasing, through Bonferroni's inequality we get
\beam \label{eq.AKP.3.39}  \notag
\PP\left[{\bf X}^{(1)}+{\bf X}^{(2)} \in x\,A\right] &\geq& \PP\left[\bigcup_{i=1}^2\left\{{\bf X}^{(i)} \in x\,A \right\}\right]\\[2mm] \notag
&\geq& \PP\left[{\bf X}^{(1)} \in x\,A\right]+ \PP\left[{\bf X}^{(2)} \in x\,A\right] - \PP\left[{\bf X}^{(1)} \in x\,A\,,\;{\bf X}^{(2)} \in x\,A\right]\\[2mm] \notag
&=& \PP\left[Y_A^{(1)}> x\right]+ \PP\left[Y_A^{(2)} >x\right] - \PP\left[Y_A^{(1)}> x\right]\,\PP\left[Y_A^{(2)} >x\right]\\[2mm] \notag
&\sim&  \overline{F_A^{(1)}}(x)+\overline{F_A^{(2)}}(x) - o\left[\overline{F_A^{(1)}}(x)\right]\\[2mm]
&\sim& \PP\left[{\bf X}^{(1)} \in x\,A\right] + \PP\left[{\bf X}^{(2)} \in x\,A\right]\,,
\eeam
as $\xto$. Hence,
\beam \label{eq.AKP.3.40} 
\liminf_{\xto} \dfrac{\PP\left[{\bf X}^{(1)}+{\bf X}^{(2)} \in x\,A\right]}{\PP\left[{\bf X}^{(1)} \in x\,A\right]+\PP\left[{\bf X}^{(2)} \in x\,A\right]} \geq 1 >0\,.
\eeam
Thus, by relations \eqref{eq.AKP.3.38} and \eqref{eq.AKP.3.40} we obtain 
\beao
\PP\left[{\bf X}^{(1)}+{\bf X}^{(2)} \in x\,A\right] \asymp \PP\left[{\bf X}^{(1)} \in x\,A\right]+\PP\left[{\bf X}^{(2)} \in x\,A\right]\,,
\eeao 
as $\xto$, which leads to relation \eqref{eq.AKP.3.35}. Via \eqref{eq.AKP.3.32} and \eqref{eq.AKP.3.35} we find
\beao
\overline{F_A^*} (x) \asymp \overline{F_A^{(1)}}(x)+\overline{F_A^{(2)}}(x)\,,
\eeao 
as $\xto$. Therefore, by $\overline{F_A^{(1)}}+\overline{F_A^{(2)}} \in \mathcal{D}$, we find  $F_A^* \in \mathcal{D}$ from where we conclude the $F_1*F_2 \in \mathcal{D}_A$. 

\item[(ii)]
We already have known that distribution class $\mathcal{D}_A$ is closed with respect to convolution, from the previous part. By \cite[Th. 2.1]{cai:tang:2004} we obtain
\beam \label{eq.AKP.3.41} 
\PP\left[Y_A^{(1)} + Y_A^{(2)} > x\right] \sim \overline{F_A^{(1)}}(x)+\overline{F_A^{(2)}}(x)\,.
\eeam
as $\xto$, when $F_A^{(1)},\,F_A^{(2)} \in \mathcal{D} \cap \mathcal{L}$. Due to Proposition \ref{pr.AKP.2.4} and relation \eqref{eq.AKP.3.41} we find
\beam \label{eq.AKP.3.42} 
\PP\left[Y_A^* > x\right] &=& \PP\left[{\bf X}^{(1)}+{\bf X}^{(2)} \in x\,A\right] \\[2mm] \notag
&\leq& \PP\left[Y_A^{(1)} + Y_A^{(2)} > x\right] \sim \overline{F_A^{(1)}}(x)+\overline{F_A^{(2)}}(x)\sim  F_1(x\,A) +F_2(x\,A)\,,
\eeam
as $\xto$. From relation \eqref{eq.AKP.3.39} and \eqref{eq.AKP.3.42} we obtain \eqref{eq.AKP.3.36}. Next, for the closure property of class $\mathcal{L}_A$, we observe that for any $y>0$ 
\beao
1&\leq& \lim_{\xto} \dfrac{\PP\left[Y_A^* > x -y\right]}{\PP\left[Y_A^* > x \right]} \\[2mm]
&=& \lim_{\xto} \dfrac{\overline{F_A^{(1)}}(x-y)+\overline{F_A^{(2)}}(x-y)}{\overline{F_A^{(1)}}(x)+\overline{F_A^{(2)}}(x)}\\[2mm]
& \leq&  \max\left\{\lim_{\xto} \dfrac{\overline{F_A^{(1)}}(x-y)}{\overline{F_A^{(1)}}(x)}\,,\;\lim_{\xto} \dfrac{\overline{F_A^{(2)}}(x-y)}{\overline{F_A^{(2)}}(x)} \right\} =1\,, 
\eeao
where in the second step we used relation \eqref{eq.AKP.3.36}. Hence, we have $F_A^* \in \mathcal{L}$, so we find $F_1*F_2 \in \mathcal{L}_A$, from where, in combination with statement $(i)$, we conclude $F_1*F_2 \in (\mathcal{D} \cap \mathcal{L})_A$.

\item[(iii)]
At first, from the second part of this theorem and the fact that $\mathcal{C}_A \subsetneq (\mathcal{D} \cap \mathcal{L})_A$, it follows that  \eqref{eq.AKP.3.36} is true. To show that $F_1*F_2 \in \mathcal{C}_A$, it is enough to obtain $F_A^* \in \mathcal{C}$. Indeed, we calculate
\beao
1\leq \lim_{b \uparrow 1}\limsup_{\xto} \dfrac{\PP\left[Y_A^* > b\,x\right]}{\PP\left[Y_A^* > x \right]} &=&\lim_{b \uparrow 1}\limsup_{\xto} \dfrac{\PP\left[Y_A^{(1)} > b\,x\right]+\PP\left[Y_A^{(2)} > b\,x\right]}{\PP\left[Y_A^{(1)} > x\right]+\PP\left[Y_A^{(2)} > x\right]}\\[2mm]
& \leq&  \max_{i \in \{1,\,2\}}\left\{\lim_{b \uparrow 1}\limsup_{\xto} \dfrac{\PP\left[Y_A^{(i)} > b\,x\right]}{\PP\left[Y_A^{(i)} > x\right]}\right\} =1\,, 
\eeao
where in the last step we took into account that $F_A^{(1)},\,F_A^{(2)} \in \mathcal{C}$, while in the second step we used relation \eqref{eq.AKP.3.36}. Therefore we conclude $F_A^* \in \mathcal{C}$.~\halmos
\end{enumerate}

\bre \label{rem.AKP.4.3}
The condition $F_A^{(1)}*F_A^{(2)} \in \mathcal{S}$ is not negligible. In \cite{leslie:1989} we find a counterexample, based on \cite{embrechts:goldie:1980}, where it was shown that class $\mathcal{S}$ is not closed with respect to convolution. Recently in \cite{leipus:siaulys:2020} were established necessary and sufficient conditions for the closure of subexponentiality with respect to convolution of two distributions from the class $\mathcal{L}(\gamma)$, for some $\gamma \geq 0$. Especially, \cite[Th. 1.1]{leipus:siaulys:2020} for the case $\gamma =0$, states that if $F_A^{(1)},\,F_A^{(2)} \in \mathcal{L}$, then the following statements are equivalent:
\begin{enumerate}
\item
$F_A^{(1)}*F_A^{(2)} \in \mathcal{S}$. 
\item
$F_{Y_A^{(1)} \vee Y_A^{(2)}} \in \mathcal{S}$.
\item
$p\,F_A^{(1)}+ (1-p)\,F_A^{(2)} \in \mathcal{S}$ for any (or, equivalently, for some) $p \in (0,\,1)$. 
\end{enumerate}
Moreover, each of these statements implies the asymptotic equivalence 
\beam \label{eq.AKP.3.a} 
\overline{F_A^{(1)}* F_A^{(2)}}(x) \sim  \overline{F_A^{(1)}}(x)+\overline{F_A^{(2)}}(x)\,, 
\eeam
as $\xto$, where by $F_{Y_A^{(1)} \vee Y_A^{(2)}}$ is denoted the distribution of the maximum. We should mention that if $F_A^{(1)},\,F_A^{(2)} \in \mathcal{S}$, then the points (1) - (3) and  \eqref{eq.AKP.3.a} are equivalent, see \cite[Cor. 1.1]{leipus:siaulys:2020}.
\ere 
 
Inspired by this result, we examine the closure property of the multivariate subexponentiality with respect to convolution. The following statement, provides necessary and sufficient conditions for the closure property of multivariate subexponentiality, that are finally reduced to corresponding one-dimensional ones. 

\bth \label{lem.AKP.3.2}
We assume that $A \in \mathscr{R}$ is a fixed set, and the claims  ${\bf X}^{(1)},\,{\bf X}^{(2)}$ are independent random vectors, with distributions $F_1,\,F_2 \in \mathcal{L}_A$ respectively. Then it holds $F_A^{(1)}*F_A^{(2)}  \in  \mathcal{S}$ if and only if  $F_1*F_2 \in  \mathcal{S}_A$.
\ethe

\pr~
($\Rightarrow$). Let $F_A^{(1)}*F_A^{(2)}  \in  \mathcal{S}$. In order to establish that $F_1*F_2 \in  \mathcal{S}_A$ we should show that $F_A^* \in \mathcal{S}$.

From relation \eqref{eq.AKP.3.39} and recalling \eqref{eq.AKP.3.32}, we obtain
\beao
\PP\left[Y_A^* >x\right] \gtrsim \PP\left[Y_A^{(1)} >x\right]+\PP\left[Y_A^{(2)} >x\right]\,,
\eeao 
as $\xto$. Thus, let ${Y_A^*}'$ be an independent copy of $Y_A^*$. Further, for $i=1,\,2$, let ${Y_A^{(i)}}'$ be an independent copy of $Y_A^{(i)}$, (with ${Y_A^{(1)}}'$ independent of $Y_A^{(2)'}$). Then we obtain
\beam \label{eq.AKP.3.21} \notag
\limsup_{\xto} \dfrac{\PP\left[Y_A^{*} + {Y_A^*}'> x\right] }{\PP\left[Y_A^{*} > x\right]} &\leq& \limsup_{\xto} \dfrac{\PP\left[Y_A^{*} + {Y_A^*}' > x\right] }{ \PP\left[Y_A^{(1)} > x\right]+ \PP\left[Y_A^{(2)} > x\right]}  \\[2mm] \notag
&\leq&\limsup_{\xto} \dfrac{\PP\left[ \left(Y_A^{(1)} + Y_A^{(2)}\right)+ \left({Y_A^{(1)}}' + {Y_A^{(2)}}'\right)> x\right] }{ \PP\left[Y_A^{(1)} > x\right]+ \PP\left[Y_A^{(2)} > x\right]} \\[2mm] \notag
&=& \limsup_{\xto} \dfrac{2\,\PP\left[ Y_A^{(1)} + Y_A^{(2)}> x\right] }{ \PP\left[Y_A^{(1)} > x\right]+ \PP\left[Y_A^{(2)} > x\right]}\\[2mm]
&=& \limsup_{\xto} \dfrac{2\,\left(\PP\left[ Y_A^{(1)} > x\right]+ \PP\left[ Y_A^{(2)} > x\right]\right)}{ \PP\left[Y_A^{(1)} > x\right]+ \PP\left[Y_A^{(2)} > x\right]}=2\,,
\eeam
where in the second step we applied an argument similar to the proof of Proposition  \ref{pr.AKP.2.4}. The third step uses the assumption $F_A^{(1)}*F_A^{(2)}  \in  \mathcal{S}$, which, by Remark \ref{rem.AKP.4.3}, implies the max-sum equivalence of $F_A^{(1)}$ and $F_A^{(2)}$.

From the other side, for non-negative and independent  random variables $Y_A^{*} $, ${Y_A^*}'$, we find
\beam \label{eq.AKP.3.22}
\liminf_{\xto} \dfrac{\PP\left[Y_A^{*} + {Y_A^*}' > x\right]}{\PP \left[Y_A^{*} > x\right]} \geq 2\,.
\eeam
Hence, from relations \eqref{eq.AKP.3.21} and \eqref{eq.AKP.3.22} we have $F_A^* \in \mathcal{S}$.

($\Leftarrow$). Let assume now $F_1*F_2 \in  \mathcal{S}_A$, that implies $F_A^* \in \mathcal{S}$. From the fact that $F_1,\,F_2 \in \mathcal{L}_A$, by \cite[Th. 3(b)]{embrechts:goldie:1980} we obtain $F_A^{(1)}*F_A^{(2)}  \in  \mathcal{L}$.

Furthermore, by Proposition \ref{pr.AKP.2.4} we find the lower bound
\beao
\overline{F_A^{(1)}*F_A^{(2)}}(x) \geq \PP\left[{\bf X}^{(1)}+{\bf X}^{(2)} \in x\,A\right]=\bF_A^*(x)\,,
\eeao
and from the other side,  since $F_A^* \in \mathcal{S}$, we can write
\beao
\overline{F_A^{(1)}*F_A^{(2)}}(x)=\PP\left[Y_A^{(1)} + Y_A^{(2)} > x\right] \leq \PP\left[Y_A^{*}+{Y_A^*}' > x\right] \sim 2\,\bF_A^*(x)\,,
\eeao
as $\xto$, where in the second step we employ the definitions of $Y_A^{(1)},\, Y_A^{(2)},\,Y_A^{*},\,{Y_A^*}'$ and the fact that the random vectors ${\bf X}^{(1)},\,{\bf X}^{(2)}$ are non-negative.

Therefore, by the last two inequalities, follows that $\overline{F_A^{(1)}*F_A^{(2)}}(x) \asymp \bF_A^*(x)$, as $\xto$, but since $F_A^* \in \mathcal{S}$ and $F_A^{(1)}*F_A^{(2)}  \in  \mathcal{L}$, from \cite[Th. 2.1(a)]{kluppelberg:1988}; see also in \cite[Prop. 3.13(ii)]{leipus:siaulys:konstantinides:2023}, we get $F_A^{(1)}*F_A^{(2)}  \in  \mathcal{S}$.
~\halmos

From Theorem \ref{lem.AKP.3.2}, since class $\mathcal{S}$ is not closed, with respect to convolution; see \cite{leslie:1989}, is implied the following result.

\bpr \label{pr.AKP.3.B}
Let  $A \in \mathscr{R}$ be a fixed set. Class $\mathcal{S}_A$ is not closed, with respect to convolution, namely  we can find two independent random vectors with distribution from class $\mathcal{S}_A$, whose convolution does NOT belong to class $\mathcal{S}_A$.
\epr

By Theorem \ref{lem.AKP.3.2}, we also get easily the following corollary, that shows an easy way to control the closure property of $\mathcal{S}_A$ with respect to convolution. 

\bco \label{th.AKP.3.3}
We assume that $A \in \mathscr{R}$ is a fixed set, and ${\bf X}^{(1)},\,{\bf X}^{(2)}$ are independent random vectors, with distributions $F_1$ and $F_2$ respectively. If $F_{1},\,F_{2}  \in  \mathcal{L}_A$, then the statements (1) - (3)  of Remark \eqref{rem.AKP.4.3} are equivalent to
$F_1*F_2 \in  \mathcal{S}_A$. Furthermore, in such a case (namely, $F_1*F_2 \in  \mathcal{S}_A$)  \eqref{eq.AKP.3.36} holds. 
\eco 

\pr~
The parts (1) - (3) of Remark \ref{rem.AKP.4.3} are equivalent, see \cite[Th. 1.1]{leipus:siaulys:2020}, and from Theorem \ref{lem.AKP.3.2} the part (1) is also equivalent
 to $F_1*F_2 \in \mathcal{S}_A$. Further, from part (1) of Remark \ref{rem.AKP.4.3} we obtain the asymptotic equivalence $\overline{F_A^{(1)}*F_A^{(2)}}(x) \sim \overline{F_A^{(1)}}(x)+ \overline{F_A^{(2)}}(x)$, as $\xto$, from where through relations \eqref{eq.AKP.3.39} and \eqref{eq.AKP.3.42} we conclude \eqref{eq.AKP.3.36}. 
~\halmos

\section{Sums of random vectors}  \label{sec.KP.4}

\subsection{Finite sums.}

Now, we study the asymptotic behavior of the entrance probability of finite sums of random vectors in the 'rare set' $x\,A$. In one-dimensional case a desired property for the subexponential distributions takes the form
\beam \label{eq.AKP.4.43}  
\PP\left[\sum_{i=1}^n X_{i} >  x \right] \sim \sum_{i=1}^n\PP\left[ X_{i} >  x \right]\,,
\eeam
as $\xto$, with fixed $n \in \bbn$. It is clear that if the random variables $X_1,\,\ldots,\,X_n$ are independent and identically distributed with common subexponential distribution, then relation \eqref{eq.AKP.4.43} holds. There exists a series of papers, where relation \eqref{eq.AKP.4.43} is studied under various conditions  for the distributions of $X_i$, $i=1,\,\ldots,\,n$, as well as for the dependence among these random variables; see \cite{ng:tang:yang:2002}, \cite{albrecher:asmussen:kotschak:2006}, \cite{chen:yuen:2009}, \cite{geluk:tang:2009} among others. It is worth to notice that in these papers there exists some 'insensitivity' with respect to dependence structure, in the sense that the heavy tails show to asymptotically override the dependence. However, there are some cases in which this 'insensitivity' with respect to dependence does not appear; see \cite{mikosch:samorodnitsky:2000a}, \cite{mikosch:samorodnitsky:2000b}, \cite{korshunov:schlegel:schmidt:2003}, \cite{foss:konstantopoulos:zachary:2007}. We keep in mind the single big jump principle, thence the multivariate form of \eqref{eq.AKP.4.43} can be given as
\beam \label{eq.AKP.4.44}  
\PP\left[{\bf S}_n \in  x\,A \right] \sim \sum_{i=1}^n\PP\left[ {\bf X}^{(i)} \in  x\,A  \right]\,,
\eeam
as $\xto$, with fixed $n \in \bbn$, where ${\bf S}_n:= {\bf X}^{(1)}+\cdots + {\bf X}^{(n)}$. Relation \eqref{eq.AKP.4.44} represents the multivariate linear single big jump principle, where the word 'linear' stems from the insensitivity of this property with respect to dimension. See also \cite[Sec. 5]{konstantinides:passalidis:2023} for more discussions in relation to the 'non-linear' case.

Usually, the study of asymptotic behavior of ${\bf S}_n$, for heavy-tailed summands, is restricted in $MRV$ class; see for example \cite{hult:lindskog:2006b}, \cite{resnick:2007}, \cite{das:fasenhartmann:2023}. In somehow bigger classes, we find in \cite{haegele:2020} and \cite{haegele:lehtomaa:2021} the study of the asymptotic
behavior of the random walk and large deviations of independent and identically distributed random vectors. In \cite[Cor. 4.10]{samorodnitsky:sun:2016} it was shown that for independent and identically distributed random vectors whose distribution belongs to $\mathcal{S}_A$, relation \eqref{eq.AKP.4.44} holds. We extend this result, by relaxing the conditions of independence and identical distribution.

Let us present three dependence structures among the random vectors, inspired by \cite{lehmann:1966}, \cite{geluk:tang:2009} and \cite{chen:yuen:2009} respectively.

\begin{assumption}[Lehmann] \label{ass.AKP.4.2}
Let $A \in \mathscr{R}$ be a fixed set. We say that the ${\bf X}^{(1)},\,\ldots,\,{\bf X}^{(n)}$, with distributions $F_1,\,\ldots,\,F_n$ respectively, are regression dependent on $A$, symbolically we write ${\bf X}^{(1)},\,\ldots,\,{\bf X}^{(n)} \in RD_A$, if $Y_{A}^{(1)},\,\ldots,\,Y_{A}^{(n)}$ are regression dependent (say $RD$), that means there exist constants $x_0$ and $C$, such that it holds
\beam \label{eq.AKP.4.46} 
\PP\left[Y_{A}^{(i)}> x_i\;\Big|\;Y_{A}^{(j)}=x_j\,, \;\text{with} \;j\in J_i\right]\leq C\,\PP\left[Y_{A}^{(i)} >x_i \right]\,,
\eeam
for $i=1,\,\ldots,\,n$, $j\in J_i$, with $\emptyset \neq J_i \subseteq \{1,\,\ldots,\,n\} \setminus \{i\}$, for $x_i \wedge x_j > x_0$, for all $j \in J_i$.
\end{assumption}

\begin{assumption}[Geluk - Tang] \label{ass.AKP.4.1}
Let $A \in \mathscr{R}$ be a fixed set. We say, that the random vectors ${\bf X}^{(1)},\,\ldots,\,{\bf X}^{(n)}$, with distributions $F_1,\,\ldots,\,F_n$, are tail asymptotically independent on $A$, symbolically ${\bf X}^{(1)},\,\ldots,\,{\bf X}^{(n)} \in TAI_A$, if $Y_{A}^{(1)},\,\ldots,\,Y_{A}^{(n)}$ are $TAI$, that means it holds
\beam \label{eq.AKP.4.45} 
\lim_{x_i\wedge x_j \to \infty } \PP \left[ Y_{A}^{(i)}> x_i\;\Big|\; Y_{A}^{(j)} > x_j \right]=0\,,
\eeam
for all $i \neq j$ with $i,\,j \in\{1,\,\ldots,\,n\}$.
\end{assumption}

\begin{assumption}[Chen - Yuen] \label{ass.AKP.4.3}
Let $A \in \mathscr{R}$ be a fixed set. We say, that the random vectors ${\bf X}^{(1)},\,\ldots,\,{\bf X}^{(n)}$, with distributions $F_1,\,\ldots,\,F_n$, are quasi asymptotically independent on $A$, symbolically ${\bf X}^{(1)},\,\ldots,\,{\bf X}^{(n)} \in QAI_A$, if the $Y_{A}^{(1)},\,\ldots,\,Y_{A}^{(n)}$ are $QAI$, namely it holds
\beam \label{eq.AKP.4.a} 
\lim_{\xto} \dfrac{\PP\left[ Y_{A}^{(i)}> x\,,\;Y_{A}^{(j)} > x \right]}{\PP\left[ Y_{A}^{(i)}> x \right]+\PP\left[ Y_{A}^{(j)} > x \right]}=0\,,
\eeam
for all $i \neq j$ with $i,\,j \in\{1,\,\ldots,\,n\}$.
\end{assumption}

Let us remind, that each random vector has arbitrarily dependent components.

\bre \label{rem.AKP.4.1}
Let us notice, that from relations \eqref{eq.AKP.4.46}, \eqref{eq.AKP.4.45}, \eqref{eq.AKP.4.a}, we can easily verify that $RD_A \subsetneq TAI_A \subsetneq QAI_A$, for any  $A \in \mathscr{R}$, and furthermore if the ${\bf X}^{(1)},\,\ldots,\,{\bf X}^{(n)}$ are independent, that means the $Y_A^{(1)},\,\ldots,\,Y_A^{(n)}$ are also independent, then they satisfy the dependence $RD_A$. Dependencies $RD$, $TAI$ and $QAI$ have been used in several applications of risk theory and risk management; see for example \cite{wang:2011}, \cite{li:2013}, \cite{chen:liu:2022}.
\ere

Now, we provide the main result.

\bth \label{th.AKP.4.1} 
\begin{enumerate}
\item[(i)]
Under Assumption \ref{ass.AKP.4.3} and if holds $F_1,\,\ldots,\,F_n \in \mathcal{C}_A$, then  the asymptotic relation \eqref{eq.AKP.4.44} is valid.

\item[(ii)]
Under Assumption \ref{ass.AKP.4.1} and if holds $F_1,\,\ldots,\,F_n \in (\mathcal{D} \cap \mathcal{L})_A$, then  the asymptotic relation \eqref{eq.AKP.4.44} is valid. 

\item[(iii)]
Under Assumption \ref{ass.AKP.4.2} and if holds $F_1,\,\ldots,\,F_n \in \mathcal{S}_A$ and $F_i*F_j \in \mathcal{S}_A$  for all $i,\,j \in \{1,\,\ldots,\,n\}$ with $i\neq j$, then  the asymptotic relation \eqref{eq.AKP.4.44} is valid. 
\end{enumerate}
\ethe

\pr~
\begin{enumerate}
\item[(i)]
By Proposition \ref{pr.AKP.2.4} and \cite[Th. 3.1]{chen:yuen:2009} we obtain
\beam \label{eq.AKP.4.47b}
\PP\left[{\bf S}_n \in  x\,A \right] \leq \PP\left[\sum_{i=1}^n Y_{A}^{(i)}> x\right]\sim \sum_{i=1}^n \PP\left[Y_{A}^{(i)}> x\right] = \sum_{i=1}^n \PP\left[{\bf X}^{(i)} \in  x\,A  \right]\,,
\eeam
as $\xto$, which provides the upper bound for \eqref{eq.AKP.4.44}. From the other side, for the lower bound, since the ${\bf X}^{(1)},\,\ldots,\,{\bf X}^{(n)}$ are non-negative random vectors and $A$ is an increasing set, using Bonferroni's inequality we find
\beam \label{eq.AKP.4.48b} \notag
\PP\left[{\bf S}_n \in  x\,A \right] &\geq& \PP\left[\bigcup_{i=1}^n \left\{ {\bf X}^{(i)} \in  x\,A \right\}\right]\\[2mm] \notag
&\geq& \sum_{i=1}^n \PP\left[ {\bf X}^{(i)} \in  x\,A \right] - \sum_{1\leq i<j\leq n} \PP\left[{\bf X}^{(i)} \in  x\,A\,,\; {\bf X}^{(j)} \in  x\,A\right]\\[2mm] \notag
&=&\sum_{i=1}^n \PP\left[Y_{A}^{(i)}> x\right]- \sum_{1\leq i<j\leq n} \PP\left[Y_{A}^{(i)}> x\,,\; Y_{A}^{(j)}> x \right]\\[2mm]
&\sim& \sum_{i=1}^n \PP\left[Y_{A}^{(i)}> x\right] = \sum_{i=1}^n \PP\left[{\bf X}^{(i)} \in  x\,A\right]\,,
\eeam
as $\xto$, where in the pre-last step due to Assumption \ref{ass.AKP.4.3}, for term of the second sum we can write 
\beao
\PP\left[Y_{A}^{(i)}> x\,,\; Y_{A}^{(j)}> x \right]= o\left(\PP\left[ Y_{A}^{(i)}> x \right]+\PP\left[ Y_{A}^{(j)}> x \right] \right)\,,
\eeao 
as $\xto$, for $i,\,j=1,\,\ldots,\,n$ and $i\neq j$. From relations \eqref{eq.AKP.4.47b} and \eqref{eq.AKP.4.48b} we find \eqref{eq.AKP.4.44}. 
\item[(ii)]
By Proposition \ref{pr.AKP.2.4} and \cite[Th. 3.1]{geluk:tang:2009} we obtain \eqref{eq.AKP.4.47b}. For the lower bound, because of the inclusion $TAI_A \subsetneq QAI_A$, the relation  \eqref{eq.AKP.4.48b} remains true. From relations \eqref{eq.AKP.4.47b} and \eqref{eq.AKP.4.48b} we find \eqref{eq.AKP.4.44}.  
\item[(iii)]
Because of the inclusion $RD_A \subsetneq QAI_A$, the relation  \eqref{eq.AKP.4.48b} remains true in this case. The upper bound of relation  \eqref{eq.AKP.4.44} follows similarly by \eqref{eq.AKP.4.47b}, where we used Proposition \ref{pr.AKP.2.4} and \cite[Th. 3.2]{geluk:tang:2009}. ~\halmos
\end{enumerate}

\subsection{Randomly stopped sums.}

In the one-dimensional case, several papers study the asymptotic behavior of the distribution tail of the randomly stopped sum
\beao
S_N:= \sum_{i=1}^N X_i\,,
\eeao
where $N$ represents a discrete, non-negative, non-degenerate to zero random variable. We note that  the summation over an empty index set is equal to zero by convention (and to ${\bf 0}$ in $d$-dimensions). For the case with heavy-tailed $X_i$, see \cite{pakes:2007}, \cite{watanabe:2008}, \cite{denisov:foss:korshunov:2010}, \cite{xu:foss:wang:2015}, among many others. In multidimensional set up, we are considering the probability $\PP\left[{\bf S}_N \in x\,A \right]$, where 
\beao
{\bf S}_N:=\sum_{i=1}^{N} {\bf X}^{(i)}\,,
\eeao 
where $N$ is a discrete, non-negative, non-degenerate to zero random variable, that is independent of $\left\{{\bf X}^{(1)},\,{\bf X}^{(2)},\,\ldots \right\}$. The asymptotic behavior of ${\bf S}_N$ can play crucial role in insurance industry. For example, we consider an insurer with $d$-lines of business, then over a concrete time interval the number of claims is unknown. For some relatively commonly usefull forms which can take the set $A$; see below in Remark \ref{rem.AKP.5.1}. Next, we extend \cite[Th. 3.37]{foss:korshunov:zachary:2013} in multivariate set up.

\bth \label{th.AKP.4.2}
We assume that $A \in \mathscr{R}$ is a fixed set, and $\{{\bf X}^{(1)},\,{\bf X}^{(2)},\,\ldots\}$, are non-negative, independent random vectors with common distribution $F\in \mathcal{S}_A$. Let $N$ be a discrete, non-negative, non-degenerate to zero random variable, such that $\E[N] <\infty$, and 
\beao
\E\left[(1 + \vep)^N\right]< \infty\,,
\eeao 
for some $\vep>0$. If further  $N$ is independent of $\{{\bf X}^{(1)},\,{\bf X}^{(2)},\,\ldots\}$, then it holds
\beam \label{eq.AKP.4.2.2}
\PP\left[{\bf S}_N \in x\,A \right] \sim  \E[N]\,F(x\,A)\,,
\eeam
as $\xto$. Furthermore, the distribution of ${\bf S}_N$ belongs to class $\mathcal{S}_A$.
\ethe

\pr~
Let assume that the right endpoint of $N$ is infinite. For the upper bound we obtain for $x>0$
\beam \label{eq.AKP.4.2.3} \notag
&&\PP\left[{\bf S}_N \in x\,A \right] = \sum_{n=1}^{\infty} \PP\left[{\bf S}_n \in x\,A \right]\,\PP[N =n] \\[2mm] \notag
&&= \sum_{n=1}^{\infty} \PP\left[\sup_{{\bf p}\in I_A} {\bf p}^{\top}\,\left( {\bf X}^{(1)}+\cdots+{\bf X}^{(n)}\right)> x\right]\,\PP[N =n] \\[2mm]
&&\leq \sum_{n=1}^{\infty} \PP\left[\sup_{{\bf p}\in I_A} {\bf p}^{\top}\, {\bf X}^{(1)}+\cdots+\sup_{{\bf p}\in I_A} {\bf p}^{\top}\,{\bf X}^{(n)}> x\right]\,\PP[N =n] \\[2mm]  \notag
&& =\sum_{n=1}^{\infty} \PP\left[\sum_{i=1}^n Y_{A}^{(i)}> x\right]\,\PP[N =n] =\PP\left[\sum_{i=1}^N Y_{A}^{(i)}> x\right] \sim  \E[N]\,\bF_A(x) =  \E[N]\,F(x\,A)\,,
\eeam
as $\xto$, where in the first step we exclude the term for $n=0$, since ${\bf 0} \notin \overline{A}$. In the sixth step we can use \cite[Th. 3.37]{foss:korshunov:zachary:2013} because of condition on $N$ and the fact that $\{ Y_{A}^{(i)}\}$ is sequence of independent, identically distributed random variables, following subexponential distribution.

From the other side, for the lower bound for \eqref{eq.AKP.4.2.2}, from \cite[Cor. 4.10]{samorodnitsky:sun:2016}; or by Theorem \ref{th.AKP.4.1} (iii), via Fatou's lemma, we obtain
\beam \label{eq.AKP.4.2.4}\notag
\PP\left[{\bf S}_N \in x\,A \right] &=& \sum_{n=1}^{\infty} \PP\left[{\bf S}_n \in x\,A \right]\,\PP[N =n] \\[2mm]
&\gtrsim& \sum_{n=1}^{\infty} n\PP\left[{\bf X} \in x\,A\right]\,\PP[N =n] = \E[N]\,F(x\,A)\,,
\eeam
as $\xto$. From relations \eqref{eq.AKP.4.2.3} and \eqref{eq.AKP.4.2.4} we get \eqref{eq.AKP.4.2.2}.

If the support of $N$ is bounded from above, with right endpoint $M$, the  relations \eqref{eq.AKP.4.2.3} and \eqref{eq.AKP.4.2.4} hold similarly with summation from one to $M$.

The membership of the distribution of ${\bf S}_N$ to class $\mathcal{S}_A$ follows by $F\in \mathcal{S}_A$, relation \eqref{eq.AKP.4.2.2}, and \cite[Prop. 4.12(a)]{samorodnitsky:sun:2016}. 
~\halmos

In the next theorem is established relation \eqref{eq.AKP.4.2.2} in the case of not necessarily independent $\{{\bf X}^{(1)},\,{\bf X}^{(2)},\,\ldots \}$, but they follow a common distribution from some classes smaller than $\mathcal{S}_A$. Additionally, we should require an alternative moment condition on the discrete random variable $N$.
  
\bth \label{th.AKP.4.3}
Let  $A \in \mathscr{R}$ be a fixed set, and $\{ {\bf X}^{(1)},\,{\bf X}^{(2)},\,\ldots\}$ be a sequence of non-negative random vectors, with common distribution $F$. Let $N$ be a non-negative, non-degenerate to zero, discrete random variable, independent of $\{{\bf X}^{(1)},\,{\bf X}^{(2)},\,\ldots\}$.
\begin{enumerate}
\item[(i)]
If the  $\{{\bf X}^{(1)},\,{\bf X}^{(2)},\,\ldots \}$ are $TAI_A$, with $F \in (\mathcal{D} \cap \mathcal{L})_A$ and it holds 
\beam \label{eq.AKP.4.3.1}
\E\left[N^{p+1} \right] < \infty \,,
\eeam
for some $p > J_{F_A}^+$, then relation \eqref{eq.AKP.4.2.2} is true, and further the distribution of ${\bf S}_N$ belongs to class $(\mathcal{D} \cap \mathcal{L})_A$.
\item[(ii)]
If the  $\{{\bf X}^{(1)},\,{\bf X}^{(2)},\,\ldots \}$ are $QAI_A$, with $F \in \mathcal{C}_A$ and \eqref{eq.AKP.4.3.1} holds, for some $p > J_{F_A}^+$, then relation \eqref{eq.AKP.4.2.2} is true, and further the distribution of ${\bf S}_N$ belongs to class $\mathcal{C}_A$.
\item[(iii)]
Under the conditions of part $(ii)$, with the restriction $F \in MRV(\alpha,\,V,\,\mu)$,
 with $\alpha \in (0,\,\infty)$, then it holds
\beam \label{eq.AKP.4.3.2}
\PP\left[{\bf S}_N \in x\,A \right] \sim \E[N]\,\mu(A)\,\bV(x) \,,
\eeam
as $\xto$.
\end{enumerate}     
\ethe

\pr~
\begin{enumerate}
\item[(i)]
Let consider $x>0$ and $M\in \bbn$. Then it follows
\beam \label{eq.AKP.4.3.3}
\PP\left[{\bf S}_N \in x\,A \right] =\left(\sum_{n=1}^M + \sum_{n=M+1}^{\infty} \right)\PP\left[{\bf S}_n \in x\,A \right] \PP[N=n]=:I_1(x,\,M) + I_2(x,\,M)\,.
\eeam
For the first term $I_1(x,\,M)$, from Theorem \ref{th.AKP.4.1}(ii) we find
\beam \label{eq.AKP.4.3.4}
I_1(x,\,M) \sim \sum_{n=1}^M \sum_{i=1}^n \PP\left[{\bf X}^{(i)} \in x\,A \right]\,\PP[N=n]= \E\left[N\,{\bf 1}_{\{N \leq M\}} \right]\,\PP\left[{\bf X} \in x\,A \right] \,.
\eeam

For the second term $I_2(x,\,M)$, choose a constant $K>0$, satisfying the inequality
\beao
\left\lfloor \dfrac x{K}\right\rfloor > M+1\,.
\eeao
Then we obtain
\beam \label{eq.AKP.4.3.5} \notag
&&I_2(x,\,M) \leq\left( \sum_{n=M+1}^{\left\lfloor x/{K}\right\rfloor} +\sum_{n=\left\lfloor x/{K}\right\rfloor +1}^{\infty}\right) \PP\left[\sum_{i=1}^n Y_A^{(i)} >x \right]\,\PP[N=n]\\[2mm]
&&=  \sum_{n=M+1}^{\left\lfloor x/{K}\right\rfloor} n\,\PP\left[Y_A >\dfrac xn \right]\,\PP[N=n] + \PP\left[N > \left\lfloor \dfrac x{K}\right\rfloor \right]\\[2mm] \notag
&&\leq C_1\,\PP\left[Y_A >x \right]\, \sum_{n=M+1}^{\left\lfloor x/{K}\right\rfloor} n^{p+1}\,\PP[N=n] + \left\lfloor \dfrac x{K}\right\rfloor ^{-(p+1)}\,\E\left[N^{p+1}\,{\bf 1}_{\{N > \left\lfloor x/{K}\right\rfloor\}} \right]\\[2mm] \notag
&&= C_1\,\PP\left[{\bf X} \in x\,A \right]\, \sum_{n=M+1}^{\left\lfloor x/{K}\right\rfloor} n^{p+1}\,\PP[N=n] + \left\lfloor \dfrac x{K}\right\rfloor ^{-(p+1)}\,\E\left[N^{p+1}\,{\bf 1}_{\{N > \left\lfloor x/{K}\right\rfloor\}} \right]\\[2mm] \notag
&&\lesssim C_1\,\PP\left[{\bf X} \in x\,A \right]\,\E\left[N^{p+1}\,{\bf 1}_{\{N > M\}} \right]\,,
\eeam
as $\xto$, where at the first step we used Proposition \ref{pr.AKP.2.4}, at the third step the  constant $C_1>0$ follows by \eqref{eq.AKP.2.55} due to the fact that $Y_A \stackrel{d}{\sim} F_A \in \mathcal{D}\cap \mathcal{L} \subsetneq \mathcal{D}$, while in the second term we applied Markov inequality, due to \eqref{eq.AKP.4.3.1}. At the last step of \eqref{eq.AKP.4.3.5}, we took into account that if $F_A \in \mathcal{D}$, then 
\beao
x^{-p} = o\left[\bF_A(x)\right]\,,
\eeao 
as $\xto$ for any $p> J_{F_A}^+$, see \cite[Lem. 3.5]{tang:tsitsiashvili:2003a}.

Hence, from \eqref{eq.AKP.4.3.4} and \eqref{eq.AKP.4.3.5}, in combination with \eqref{eq.AKP.4.3.1}, for any $\delta >0$, we can find some large enough $M \in \bbn$, such that
\beam \label{eq.AKP.4.3.6} \notag
&&I_1(x,\,M) \sim \left(\E\left[N\right]\,-\E\left[N\,{\bf 1}_{\{N > M\}}\right] \right)\,F(x\,A ) \lesssim (1+\delta)\,\E\left[N\right]\,F(x\,A )\,, \\[2mm]
&&I_2(x,\,M) \lesssim \delta\, \E\left[N\right]\,F(x\,A )\,,
\eeam
as $\xto$. From \eqref{eq.AKP.4.3.3}, \eqref{eq.AKP.4.3.6} and the arbitrariness in the choice of $\delta>0$, we obtain upper bound in \eqref{eq.AKP.4.2.2}. 

For the lower bound in \eqref{eq.AKP.4.2.2}, because $\PP({\bf S}_n \in x\,A) \gtrsim n\,F(x\,A)$; see in Theorem \ref{th.AKP.4.1} (ii) the proof of \eqref{eq.AKP.4.48b}, via Fatou's lemma, we obtain
\beao
&&\PP\left[{\bf S}_N \in x\,A \right] = \sum_{n=1}^{\infty} \PP\left[{\bf S}_n \in x\,A \right]\,\PP[N=n] \\[2mm] 
&&\gtrsim \sum_{n=1}^{\infty} n\,\PP\left[{\bf X} \in x\,A \right]\,\PP[N=n]\sim \E[N]\,\PP\left[ {\bf X}\in x\,A \right] \,.
\eeao 
where at the second step we took into consideration that the ${\bf X}^{(i)}$, ${\bf X}^{(j)}$ are $TAI_A$. Hence, \eqref{eq.AKP.4.2.2} is true.

Further, since $F \in (\mathcal{D}\cap \mathcal{L})_A$, and due to \eqref{eq.AKP.4.2.2}, the distribution of the ${\bf S}_N$ belongs to $(\mathcal{D}\cap \mathcal{L})_A$, as it is implied by Proposition \ref{pr.AKP.2.2}(iv).
\item[(ii)]
It follows by similar methodology, using Theorem \ref{th.AKP.4.1}(i) for $I_1(x,\,M)$, and Proposition \ref{pr.AKP.2.2}(ii), to prove that the distribution of the ${\bf S}_N$ belongs to $\mathcal{C}_A$.
\item[(iii)]
\eqref{eq.AKP.4.2.2} is true, due to part $(ii)$. This, in combination with \eqref{eq.AKP.2.6} leads to \eqref{eq.AKP.4.3.2}, since $\mu(A) \in (0,\,\infty)$, for all  $A \in \mathscr{R}$. ~\halmos
\end{enumerate}

\section{Asymptotic behavior of aggregate claims.}

Now, we give an application in risk theory. More precisely, we want to find an asymptotic expression for the probability of the entrance of discounted aggregate claims  into some rare set $x\,A$, for a multivariate Poisson risk model, with claim distributions from the class $\mathcal{S}_A$ and permitting for the insurer, risk-free and risky investments simultaneously, through a general stochastic process.

In risk theory  there exist several papers considering the asymptotic behavior of discounted aggregate claims; or even the ruin probability, with claims from the class of subexponential distributions in one-dimensional or two-dimensional risk models; see for example \cite{asmussen:1998}, \cite{li:tang:wu:2010}, \cite{jiang:wang:chen:xu:2015}, \cite{yuan:lu:2023}. However in most multivariate models we meet as basic assumption the $MRV$ class for the distribution of claims, and additionally the asymptotic dependence among the components of each claim vector; see for example \cite{konstantinides:li:2016} for risk model with constant interest force and \cite{li:2016}, \cite{cheng:konstantinides:wang:2022}, \cite{cheng:konstantinides:wang:2024} for risk models with both risk-free and risky investments. In the case of larger claim distribution classes, we refer \cite{samorodnitsky:sun:2016}, where a renewal risk model without interest force and with multivariate subexponential claims is examined.

Motivated by the modern insurance industry, having in mind that most of the insurance companies invest their surpluses, we consider a multivariate risk model, which covers a general price process of the financial portfolio. The following concept is inspired by \cite{chen:liu:2023}, and generalizes in some sense their result in multidimensional set up.

We consider an insurer with $d$ lines of business, with claims ${\bf X}^{(1)},\,{\bf X}^{(2)},\,\ldots$, that are independent, identically distributed, non-negative random vectors, that arrive at moments $0 < \tau_1 < \tau_2 < \cdots$, which constitute a Poisson counting process $N_t:=\sup\{ n\in \bbn\;:\;\tau_n \leq t\}$, $\sup \emptyset = 0$, with intensity $\lambda > 0$.  The price process of the investment portfolio has the form $\{e^{R_t}\,,\;t\geq 0 \}$, where by $\{R_t\,,\;t\geq 0\}$ we denote a general real-valued stochastic process, which represents the logarithmic return of the prices of the investment portfolio. Hence, the insurer's discounted aggregate claims, up to moment $t>0$, are provided by the relation
\beam \label{eq.AKP.5.50}
{\bf D}(t) = \sum_{i=1}^{N_t}{\bf X}^{(i)} e^{-R_{\tau_i}}=\left( 
\sum_{i=1}^{N_t} X_{1}^{(i)}\,e^{-R_{\tau_{i}} }\,,\;
\cdots \,,\; \sum_{i=1}^{N_t} X_{d}^{(i)}\,e^{-R_{\tau_{i}}} \right)^{\top}\,.
\eeam 

Before reaching to the main result we present some preliminary lemmas. The first two seem to have value on their own merit. In the first we verify that the presence of the multivariate linear single big jump remains invariable under some scale mixtures, that are upper bounded and non-negative. Let us denote the scale mixture sum as
\beam \label{eq.AKP.5.54}
{\bf S}_n^{\Theta} = \sum_{i=1}^{n} \Theta^{(i)}\,{\bf X}^{(i)} \,.
\eeam
Recall that, as analogue to \eqref{eq.AKP.3.24}, we define 
\beao
M_A^{(i)}:=\sup\left\{u\;:\;\Theta^{(i)}\,{\bf X}^{(i)} \in u\,A \right\}\,,
\eeao 
for $i=1,\,\ldots,\,n$.

\ble \label{lem.AKP.5.1}
We assume that $A \in \mathscr{R}$ is a fixed set, and the claims  ${\bf X}^{(1)},\,\ldots,\,{\bf X}^{(n)}$ are independent, identically distributed random vectors, with common distribution $F \in \mathcal{S}_A$. Let $\Theta^{(1)},\,\ldots,\,\Theta^{(n)}$ be arbitrarily dependent, non-negative and non-degenerate to zero random variables, which are uniformly upper bounded by $0 < b <\infty$. If $\Theta^{(1)},\,\ldots,\,\Theta^{(n)}$ are independent of ${\bf X}^{(1)},\,\ldots,\,{\bf X}^{(n)}$, then
\beam \label{eq.AKP.5.55}
\PP\left[ {\bf S}_n^{\Theta} \in x\,A \right] \sim \sum_{i=1}^n \PP\left[ \Theta^{(i)} \,{\bf X}^{(i)} \in x\, A\right]\,, 
\eeam
as $\xto$.
\ele

\pr~
On the one hand, for the upper bound of \eqref{eq.AKP.5.55} we obtain the following
\beao
&&\PP\left[ {\bf S}_n^{\Theta} \in x\,A \right]= \PP\left[\sum_{i=1}^n  \Theta^{(i)} \,{\bf X}^{(i)} \in x\, A\right] =\PP\left[ \sup_{{\bf p}\in I_A}{\bf p}^{\top} \,\left( \sum_{i=1}^n  \Theta^{(i)} \,{\bf X}^{(i)} \right) > x\right] \\[2mm]
&&\leq \PP\left[ \sup_{{\bf p}\in I_A}{\bf p}^{\top} \,\Theta^{(1)} \,{\bf X}^{(1)} +\cdots+ \sup_{{\bf p}\in I_A}{\bf p}^{\top} \,\Theta^{(n)} \,{\bf X}^{(n)} > x\right]= \PP\left[ \sum_{i=1}^n \Theta^{(i)} \,Y_A^{(i)} >x \right] \\[2mm]
&& \sim \sum_{i=1}^n \PP\left[ \Theta^{(i)} \,Y_A^{(i)} >x \right]= \sum_{i=1}^n \PP\left[ M_A^{(i)} >x \right]= \sum_{i=1}^n  \PP\left[ \Theta^{(i)} \,{\bf X}^{(i)} \in x\, A\right]\,,
\eeao
as $\xto$, where in the fifth step we used \cite[Th. 1]{tang:yuan:2014}, while in the pre-last step we took into consideration the relation
\beam \label{eq.AKP.5.56}  \notag
\PP\left[ M_A^{(i)} >x \right]&=&\PP\left[ \sup_{{\bf p}\in I_A}{\bf p}^{\top} \,\left(  \Theta^{(i)} \,{\bf X}^{(i)} \right) > x\right]\\[2mm]
&=&\PP\left[ \Theta^{(i)} \,\sup_{{\bf p}\in I_A}{\bf p}^{\top} {\bf X}^{(i)} > x \right]=\PP\left[  \Theta^{(i)} \,Y_A^{(i)} >x \right]\,, 
\eeam
for $i=1,\,\ldots,\,n$. 

On the other hand, for the lower bound of \eqref{eq.AKP.5.55}, since $A$ is an increasing set, and $ \Theta^{(i)}$, ${\bf X}^{(i)}$ are non-negative, through Bonferroni's inequality we obtain
\beao
\PP\left[ {\bf S}_n^{\Theta} \in x\,A \right] &\geq& \PP\left[\bigcup_{i=1}^n \left\{  \Theta^{(i)} \,{\bf X}^{(i)} \in x\, A \right\} \right]\\[2mm]
&\geq& \sum_{i=1}^n  \PP\left[ \Theta^{(i)} \,{\bf X}^{(i)} \in x\, A \right] -\sum_{1\leq i<j\leq n}  \PP\left[ \Theta^{(i)} \,{\bf X}^{(i)} \in x\, A\,,\;  \Theta^{(j)} \,{\bf X}^{(j)} \in x\, A \right]\\[2mm]
& \geq&  \sum_{i=1}^n  \PP\left[ M_A^{(i)} >x \right] -\sum_{1\leq i<j\leq n} \PP\left[ M_A^{(i)} >x \right]\,\PP\left[  Y_A^{(j)} > \dfrac xb \right]\\[2mm]
& \sim& \left[1-o(1)\right]\,\sum_{i=1}^n  \PP\left[ M_A^{(i)} >x \right] \sim \sum_{i=1}^n \PP\left[ \Theta^{(i)} \,{\bf X}^{(i)} \in x\, A \right]\,,
\eeao 
as $\xto$.
~\halmos

The next lemma establishes a Kesten-type bound for the scale mixture sum of \eqref{eq.AKP.5.54}, thereby extending \cite[Prop. 4.12 (c)]{samorodnitsky:sun:2016}. 

\ble \label{lem.AKP.5.2}
Under the conditions of Lemma \ref{lem.AKP.5.1}, for any $\vep>0$, there exists constant $K=K(\vep)>0$, such that it holds
\beao
\PP\left[ {\bf S}_n^{\Theta} \in x\,A \right] \leq K\,(1+\vep)^n\, \sum_{i=1}^n \PP\left[ \Theta^{(i)} \,{\bf X}^{(i)} \in x\,A\right]\,, 
\eeao
for all $n \in \bbn$ and all $x \geq 0$.
\ele

\pr~
For any $x>0$, we obtain the following inequalities
\beao
\PP\left[ {\bf S}_n^{\Theta} \in x\,A \right] &=&\PP\left[ \sup_{{\bf p}\in I_A}{\bf p}^{\top} \,\left( \Theta^{(1)} \,{\bf X}^{(1)}+\cdots+ \Theta^{(n)} \,{\bf X}^{(n)} \right) > x\right]\\[2mm]
&\leq& \PP\left[ \sup_{{\bf p}\in I_A}{\bf p}^{\top} \,\Theta^{(1)} \,{\bf X}^{(1)}+\cdots+ \sup_{{\bf p}\in I_A}{\bf p}^{\top} \,\Theta^{(n)} \,{\bf X}^{(n)} > x\right]\\[2mm]
&=& \PP\left[ \sum_{i=1}^n \Theta^{(i)} \,Y_A^{(i)} >x \right] \leq K\,(1+\vep)^n\, \sum_{i=1}^n \PP\left[ \Theta^{(i)} \,Y_A^{(i)} > x\right]\\[2mm]
&=&  K\,(1+\vep)^n\, \sum_{i=1}^n \PP\left[M_A^{(i)} >x \right] = K\,(1+\vep)^n\, \sum_{i=1}^n \PP\left[ \Theta^{(i)} \,{\bf X}^{(i)} \in x\, A\right]\,, 
\eeao
where in the forth step we have applied \cite[Th. 1.2]{chen:2020} and in the penultimate step we have used relation \eqref{eq.AKP.5.56}.
~\halmos

Now, we remind a classical result on Poisson process; see for example \cite[Th. 5.2]{ross:2014}.

\ble \label{lem.AKP.5.3}
We assume that $0< \tau_1< \tau_2 < \ldots$ are jump epochs of a Poisson process $\{N_t\,,\;t \geq 0\}$. Let $U_1,\,U_2,\,\ldots$ be a sequence of independent, identically distributed random variables, with uniform distribution over the interval $(0,\,1)$, that are independent of  $\{N_t\,,\; t \geq 0\}$. For any arbitrarily constant $T>0$ and any $n \in \bbn$ it holds $\left\{ (\tau_1,\,\ldots,\,\tau_n)\;|\;N_T=n\right\} \stackrel{d}{=} \left(T\,U_{(1,n)},\,\ldots,\,T\,U_{(n,n)}\right)$, where the symbol $\stackrel{d}{=}$ means equality in distribution and the sequence $U_{(1,n)},\,\ldots,\,U_{(n,n)}$ represents the order statistics of $U_{1},\,\ldots,\,U_{n}$.
\ele

Let us formulate the assumptions for the risk model \eqref{eq.AKP.5.50}.

\begin{assumption} \label{ass.AKP.5.1}
For some fixed set $A \in \mathscr{R}$, let the claims  $\{{\bf X},\,{\bf X}^{(1)},\,{\bf X}^{(2)},\,\ldots \}$, be a sequence of independent, identically distributed, non-negative random vectors, with common distributions $F \in \mathcal{S}_A$. Further we assume that $\{{\bf X},\,{\bf X}^{(1)},\,{\bf X}^{(2)},\,\ldots \}$, $\{N_t\,,\;t\geq 0\}$ and $\left\{e^{R_t}\,,\;t\geq 0 \right\}$ are mutually independent.
\end{assumption}

Let us observe that the $d$ lines of business, at each jump epoch $\tau_i$, are probably dependent, that means the vector ${\bf X}^{(i)}=\left(X_{1}^{(1)},\,\ldots,\, X_{d}^{(1)}\right)^{\top}$ may has dependent components through $\mathcal{S}_A$, which contains cases of asymptotic independence and asymptotic dependence. In opposite to the most of multidimensional risk models with $MRV$ claims, where the requirement of asymptotic dependence very often is needed. Some of the components of ${\bf X}^{(i)}$ can be zero, but not all of them (in order to have renewal epoch at $\tau_i$). Finally, from the fact that the model has common counting process at the $d$ lines of business we obtain also another source of dependence between the lines.

\begin{assumption} \label{ass.AKP.5.2}
Let $T > 0$ be finite. The logarithmic price process of an investment portfolio $\left\{R_t\,,\;0< t\leq T\right\}$, is bounded away from $-\infty$, namely there exists a constant $C=C(T) < \infty$, such that satisfies the condition
\beao
\PP\left[ \inf_{0\leq t \leq T } R_t > -C \right] =1\,. 
\eeao
\end{assumption}

Assumption \ref{ass.AKP.5.2} is general enough to permit to the insurer to have both risk-free and risky investments simultaneously; see in \cite[Exam. 1]{chen:liu:2023}. The Assumption \ref{ass.AKP.5.2} is clearly dependent of the choice of $T >0$ in several cases.  

\bth \label{th.AKP.5.1}
We assume Assumption \ref{ass.AKP.5.1} and Assumption \ref{ass.AKP.5.2} are valid in risk model \eqref{eq.AKP.5.50}. Then it holds
\beam \label{eq.AKP.5.53}
\PP\left[ {\bf D}(T) \in x\,A\right] \sim \lambda\,\int_0^T \PP \left[ {\bf X}\,e^{-R_{t}} \in x\,A \right]\,dt \,,
\eeam
as $\xto$.
\ethe

\bre \label{rem.AKP.5.1}
In the one-dimensional case, relation \eqref{eq.AKP.5.53} is interesting for $A=(1,\,\infty)$, since then it represents the tail of discounted aggregate claims. However, in the multivariate risk models there are several interesting cases. For example if the set $A$ is of the form \eqref{eq.AKP.2.12}, then the \eqref{eq.AKP.5.53} express the probability that the sum of the discounted aggregate claims of the $d$-lines of business, exceeds the initial capital. As we see in \cite{samorodnitsky:sun:2016} the set $A$ can take the form $\{{\bf x}\;:\;x_i>b_i\,,\;\text{for some}\; i=1,\,\ldots,\,d\}$, where $b_1,\,\ldots,\,b_d>0$, in which we are interested in the probability that one business line to exceed its initial capital, accordingly. For a depiction of several ruin probabilities in multivariate risk models see in \cite[Sec. 1]{cheng:yu:2019}
\ere

\pr~
By the law of total probability, and the independence of the counting process $\{N_t\,,\; t \geq 0\}$ from the other sources of randomness we obtain for any $x>0$,
\beao
\PP\left[{\bf D}(T) \in x\,A\right]=\sum_{n=1}^{\infty} \PP\left[\sum_{i=1}^n {\bf X}^{(i)}\,e^{-R_{\tau_i}} \in x\,A \right]\,\PP\left[N_T=n\right]\,.
\eeao
We now define the $U_1,\,U_2,\,\ldots$ as a sequence of independent identically uniformly distributed over the interval $(0,\,1)$ random variables, which are independent of all the other sources of randomness. Then by Lemma \ref{lem.AKP.5.3} we obtain
\beam \label{eq.AKP.5.58} \notag
\PP[{\bf D}(T) \in x\,A]&=&\sum_{n=1}^{\infty} \PP\left[\sum_{i=1}^n {\bf X}^{(i)}\,e^{-R_{TU_{(i,n)}}} \in x\,A \right]\,\PP[N_T=n]\\[2mm]
&=&\sum_{n=1}^{\infty} \PP\left[\sum_{i=1}^n {\bf X}^{(i)}\,e^{-R_{TU_{i}}} \in x\,A \right]\,\PP[N_T=n]\,, 
\eeam
where in the last step we used Assumption \ref{ass.AKP.5.1}. Further the  
\beao
\sum_{i=1}^n {\bf X}^{(i)}\,e^{-R_{TU_{i}}}\,, 
\eeao
is a scale mixture sum, as in relation \eqref{eq.AKP.5.54}, with $\Theta^{(i)} =e^{R_{TU_{i}}}$, for $i=1,\,\ldots,\,n$. From Assumption  \ref{ass.AKP.5.2} we have that the  $\Theta^{(i)} $ are non-negative and upper bounded by $e^{C}$. Hence all the assumptions of Lemma \ref{lem.AKP.5.1} and Lemma \ref{lem.AKP.5.2} are satisfied, namely by the first one we find
\beam \label{eq.AKP.5.59}  
\PP\left[\sum_{i=1}^n {\bf X}^{(i)}\,e^{-R_{TU_{i}}} \in x\,A \right]
\sim \sum_{i=1}^n  \PP\left[{\bf X}^{(i)}\,e^{-R_{TU_{i}}} \in x\,A \right]=n\,\PP\left[{\bf X}\,e^{-R_{TU}} \in x\,A \right]\,, 
\eeam
as $\xto$, for $n \in \bbn$, where $U$ is a copy of $U_1,\,\ldots,\,U_n$.

Furthermore, by Lemma \ref{lem.AKP.5.2}, and for any $\vep>0$, there exists a constant $K=K(\vep)>0$ such that it holds
\beam \label{eq.AKP.5.60}  \notag 
\PP\left[\sum_{i=1}^n {\bf X}^{(i)}\,e^{-R_{TU_{i}}} \in x\,A \right]
&\leq& K\,(1+\vep)^n\, \sum_{i=1}^n \PP\left[ {\bf X}^{(i)}\,e^{-R_{TU_{i}}} \in x\,A \right]\\[2mm]
&=&K\,(1+\vep)^n\, n \PP\left[ {\bf X}\,e^{-R_{TU}} \in x\,A \right]
\,, 
\eeam
for any $n\in \bbn$ and $x\geq 0$ and additionally, since the process $\{N_t\,,\;t\geq 0\}$ is Poisson, we can find some $\vep >0$, such that it holds
\beao
\sum_{n=1}^{\infty} n\,(1+\vep)^n\,\PP\left[N_T=n\right]=\E\left[N_T\,(1+\vep)^{N_T} \right] < \infty\,, 
\eeao
for $T \in (0,\,\infty)$ (see for example \cite[Th. 1]{kocetova:leipus:siaulys:2009}). From this, together with relation \eqref{eq.AKP.5.60} we obtain
\beao
\dfrac{\PP\left[{\bf D}(T) \in x\,A\right]}{\PP\left[ {\bf X}\,e^{-R_{TU}} \in x\,A \right]} \leq K\,\E\left[N_T\,(1+\vep)^{N_T} \right] < \infty
\eeao

Therefore, we can apply the dominated convergence theorem on \eqref{eq.AKP.5.58} and through \eqref{eq.AKP.5.59} we have
\beao
\lim_{\xto}\dfrac{\PP\left[{\bf D}(T) \in x\,A\right]}{\PP\left[ {\bf X}\,e^{-R_{TU}} \in x\,A \right]} &=& \sum_{n=1}^{\infty}\lim_{\xto}\dfrac{\PP\left[\sum_{i=1}^n {\bf X}^{(i)}\,e^{-R_{TU_{i}}} \in x\,A \right]}{\PP\left[ {\bf X}\,e^{-R_{TU}} \in x\,A \right]}\,\PP\left[N_T=n\right]\\[2mm]
&=&\sum_{n=1}^{\infty} n\,\PP[N_T=n]=\E\left[N_T\right] =\lambda\,T\,,
\eeao
from where follows, through integration with respect to $U$, the asymptotic relation
\beao
\PP[{\bf D}(T) \in x\,A] \sim \lambda\,T\,\PP\left[ {\bf X}\,e^{-R_{TU}} \in x\,A \right] \sim \lambda\,\int_0^T\,\PP\left[ {\bf X}\,e^{-R_{t}} \in x\,A \right]\,dt\,,
\eeao
as $\xto$, which is the desired result.
~\halmos

\noindent \textbf{Acknowledgments.} 
We are grateful to an anonymous referee for many substantianl comments that improved the paper. We also feel the pleasant duty to express our gratitude to T. Mikosch and M. Steffensen, for their stimulating discussions, during the visit of the first author to Copenhagen.


\begin{thebibliography}{99}






\bibitem{asmussen:1998}
{\sc Asmussen, S.}\ (1998)
Subexponential asymptotics for stochastic processes: extremal behaviour, stationary distributions and first passage probabilities.
{\em Ann. Appl. Probab.}, \textbf{8}, 354--374.

\bibitem{asmussen:albrecher:2010}
{\sc Asmussen, S., Albrecher, H.}\ (2010)
{\em Ruin Probabilities.} 
World Scientific, Singapore, 2nd ed.

\bibitem{asmussen:steffensen:2020}
{\sc Asmussen, S., Steffensen, M.}\ (2020)
{\em Risk and Insurance: A graduate text} 
Springer, Cham.


\bibitem{albrecher:asmussen:kotschak:2006}
{\sc Albrecher, H., Asmussen, S., Kortschak, D.}\ (2006)
Tail asymptotics for the sums of two heavy-tailed dependent risks. 
{\em Extremes}, \textbf{9}, no. 2,  107--130.




\bibitem{basrak:davis:mikosch:2002a}
{\sc Basrak, B., Davis, R.A., Mikosch, T.}\ (2002) 
Regular variation of GARCH processes.
{\em  Stoch. Process. Appl.} \textbf{99}, no. 1, 95--115.



\bibitem{bingham:goldie:teugels:1987} 
{\sc Bingham. N.H., Goldie, C.M., Teugels, J.L.} \ (1987)
{\em Regular Variation}
Cambridge University Press, Cambridge.




\bibitem{cai:tang:2004}
{\sc Cai, J., Tang, Q.}\  (2004)
On max-sum equivalence and convolution closure of heavy-tailed distributions and their applications. 
{\em J. Appl. Probab.} \textbf{41}, 117--130.

\bibitem{chen:2020}
{\sc Chen, Y.}\  (2020)
A Kesten-type bound for sums of randomly weighted subexponential random variables. 
{\em Stat. Probab. Lett.} \textbf{158}, 108661.





\bibitem{chen:liu:2022}
{\sc Chen, Y., Liu, J.}\  (2022)
An asymptotic study of systemic expected shortfall and marginal expected shortfall. 
{\em Insur. Math. Econom.} \textbf{105}, 238--251.

\bibitem{chen:liu:2023}
{\sc Chen, Y., Liu, J.}\  (2023)
An asymptotic result on catastrophe insurance losses. 
{\em North Amer. Actur. J.} \textbf{28}, 426--437.

\bibitem{chen:liu:2024}
{\sc Chen, Y., Liu, J.}\  (2024)
Asymptotic capital allocation based on the higher moment risk measure. 
{\em Europ. Actuar. J.} \textbf{14}, 657--684.



\bibitem{chen:yuen:2009}
{\sc Chen, Y., Yuen, K.C.}\ (2009)
Sums of pairwise quasi-asymptotic independent random variables with consistent variation.
{\em Stochastic Models}, 25, 76--89.




\bibitem{cheng:yu:2019}
{\sc Cheng, D., Yu, C.}\  (2019)
Uniform asymptotics for the ruin probabilities in a bidimensional renewal risk model with strongly subexponential claims.
{\em Stochastics} \textbf{91}, Vol 1. 643--656.

\bibitem{cheng:konstantinides:wang:2022}
{\sc Cheng, M., Konstantinides, D.G., Wang, D}\ (2022)
Uniform asymptotic estimates in a time-dependent risk model with general investment returns and multivariate regularly varying claims. 
{\em Appl. Math. and Comput.}, \textbf{434}, 127436.

\bibitem{cheng:konstantinides:wang:2024}
{\sc Cheng, M., Konstantinides, D.G., Wang, D}\ (2024)
Multivariate regular varying insurance and financial risks in $d$-dimensional risk model. 
{\em J. Appl. Probab.}, \textbf{61}, no. 4, 1319 -- 1342.

\bibitem{chistyakov:1964}
{\sc Chistyakov, V.P.}\ (1964)
A theorem on sums of independent positive random variables and its applications to branching random processes.
{\em Theory Probab. Appl.}, \textbf{9}, 640--648.

\bibitem{cline:1994} 
{\sc Cline, D.B.H.}\ (1994)
Intermediate regular and $\Pi$ variation.
{\em Proc. Lond. Math. Soc.}, \textbf{68}, 594--611.

\bibitem{cline:resnick:1992} 
{\sc Cline, D.B.H., Resnick, S.}\ (1992)
Multivariate subexponential distributions.
{\em Stoch. Process. Appl.}, \textbf{42}, no.1, 49--72.

\bibitem{cline:samorodnitsky:1994} 
{\sc Cline, D.B.H., Samorodnitsky, G.}\ (1994)
Subexponentiality of the product of independent random variables.
{\em Stoch. Process. Appl.}, \textbf{49}, 75--98.


\bibitem{daley:omey:vesilo:2007} 
{\sc Daley, D.J., Omey, E., Vesilo, R.}\ (2007) 
The tail behaviour of random sums of subexponential random variables and vectors. 
{\em Extremes}, \textbf{10}, 21--39.

\bibitem{das:fasenhartmann:2023}
{\sc Das, B., Fasen-Hartmann, V.}\ (2023)
Aggregating heavy-tailed random vectors: from finite sums to Levy processes.
{\em Preprint, arXiv:2301.10423v1}.

\bibitem{denisov:foss:korshunov:2010}
{\sc Denisov, D., Foss, S., Korshunov, D.}\ (2010)
Asymptotics of randomly stopped sums in the presence of heavy tails.
{\em Bernoulli} \textbf{16}, 971--994.





\bibitem{embrechts:goldie:1980}
{\sc Embrechts, P., Goldie, C. M.}\ (1980)
On closure and factorization properties of subexponential and related distributions. 
{\em J. Austral. Math. Soc. (Ser. A)}, \textbf{29}, 243--256.


\bibitem{embrechts:klueppelberg:mikosch:1997}
{\sc Embrechts, P., Kl\"{u}pellberg, C. and Mikosch, T.}\ (1997) 
{\em  Modelling Extremal Events for Insurance and Finance.} 
Springer, New York.


\bibitem{feller:1969}
{\sc Feller, W.}\ (1969)
One-sided\quad analogues\quad of\quad  Karamata's\quad regular\quad  variation.
{\em L' enseignement Math\'{e}matique}, \textbf{15}, 107--121.


\bibitem{foss:korshunov:zachary:2013} 
{\sc Foss, S., Korshunov, D., Zachary, S.} \ (2013)
{\em An Introduction to Heavy-Tailed and Subexponential Distributions.}
Springer, New York, 2nd ed.

\bibitem{foss:konstantopoulos:zachary:2007} 
{\sc Foss, S., Konstantopoulos, T., Zachary, S.} \ (2007)
{\em Discrete and continuous time moluled random walks with heavy-tailed increments.}
{\em J. Theor. Probab.}, \textbf{20}, no. 3, 581--612.

\bibitem{fougeres:mercadier:2012}
{\sc Fougeres, A., Mercadier, C.}\ (2012)
Risk measures and multivariate extensions of Breiman’s theorem.
{\em J. Appl. Probab.},  \textbf{49}, no. 2, 364--384.



\bibitem{geluk:tang:2009}
{\sc Geluk, J., Tang, Q.}\ (2009)
Asymptotic tail probabilities of sums of dependent subexponential random variables. 
{\em J. Theor. Probab.}, \textbf{22}, 871--882.


\bibitem{goldie:1978}
{\sc Goldie, C.M.}\ (1978)
Subexponential distributions and dominated variation tails 
{\em J. Appl. Probab.}, \textbf{15}, 440--442.

\bibitem{haan:resnick:1981} 
{\sc Haan, L. de, Resnick, S.}\ (1981)
On the observation closet to the origin.
{\em  Stoch. Process. Appl.}, \textbf{11}, no. 3, 301--308.

\bibitem{haegele:2020}
{\sc Haegele, M.}\ (2020)
Precise asymptotics of ruin probabilities for a class of multivariate heavy-tailed distributions. 
{\em Stat. Prob. Lett.}, \textbf{166}, 108871.

\bibitem{haegele:lehtomaa:2021} 
{\sc Haegele, M., Lehtomaa, J.}\ (2021)
Large deviations for a class of multivariate heavy-tailed risk processes used in insurance and finance.
{\em  J. Risk Fin. Manag.}, \textbf{14}, 202.



\bibitem{hult:lindskog:2006b}
{\sc Hult, H., Lindskog, F.}\ (2006)
On regular variation for infinitely divisible random vectors and additive processes.
{\em Adv. Appl. Probab.}, \textbf{38}, 134--148.






\bibitem{jiang:wang:chen:xu:2015}
{\sc Jiang, T., Wang, Y., Chen, Y., Xu, H.}\  (2015)
Uniform asymptotic estimate for finite-time ruin probabilities of a time-dependent bidimensional renewal model. 
{\em Insur. Math. Econom.}, \textbf{64}, 45--53.



\bibitem{kaleta:ponikowski:2022} 
{\sc Kaleta, K., Ponikowski, D.} \ (2022)
On directional convolution equivalent densities.
{\em Elecrt. J. Probab.}, \textbf{27}, 1--19.

\bibitem{karamata:1933} 
{\sc Karamata, J.} \ (1933)
Sur un mode de croissance r\'{e}guli\`{e}re. Th\'{e}or\`{e}mes fondamentaux.
{\em Bulletin de la Soci\'{e}t\'{e} Math\'{e}matique de France}, \textbf{61}, 55--62.

\bibitem{kapnova:palmowski:2022} 
{\sc Kapnova, V., Palmowski, Z.} \ (2022)
Subexponential potential asymptotics with applications.
{\em Adv. Appl. Probab.}, \textbf{54}, 1--25.


\bibitem{kluppelberg:1988}
{\sc Kl\"{u}pellberg, C.}\ (1988) 
Subexponential distributions and integrated tails. 
{\em J. Appl. Probab.}, \textbf{25}, 132--141.

\bibitem{kocetova:leipus:siaulys:2009}
{\sc Ko{\v c}etova, J., Leipus, R.,  \v{S}iaulys, J.}\ (2009)
A property of the renewal counting process with application to the finite-time ruin probability.
{\em Lith. Math. J.}, \textbf{49}, no.1, 55--61.

\bibitem{konstantinides:2018} 
{\sc Konstantinides, D.G.} \ (2018)
{\em Risk Theory. A Heavy Tail Approach.}
World Scientific, New Jersey.

\bibitem{konstantinides:leipus:siaulys:2022}
{\sc Konstantinides, D.G., Leipus, R., \v{S}iaulys, J.}\ (2022) 
A note on product-convolution for generalized subexponential distributions. 
{\em Non. Anal. Mod. Contr.}, \textbf{27}, 1054--1067.


\bibitem{konstantinides:li:2016}
{\sc Konstantinides, D.G., Li, J.}\ (2016)
Asymptotic ruin probabilities for a multidimensional renewal risk model with multivariate regularly varying claims.
{\em Insur. Math. Econom.}, \textbf{69}, 38--44.

\bibitem{konstantinides:liu:passalidis:2025} 
{\sc Konstantinides, D.G., Liu, J., Passalidis, C.D.} \ (2026)
Uniform asymptotics for a multidimensional renewal risk model with multivariate subexponential claims.
{\em Scand. Actuar. J.}, p. 1 -- 21.\\ 
DOI: 10.1080/03461238.2025.2584008. 

\bibitem{konstantinides:passalidis:2024b} 
{\sc Konstantinides, D.G., Passalidis, C.D.} \ (2024a)
Closure properties and heavy tails: random vectors in the presence of dependence
{\em Preprint, arXiv:2402.09041}.


\bibitem{konstantinides:passalidis:2024c} 
{\sc Konstantinides, D.G., Passalidis, C.D.} \ (2025a)
A new approach in two-dimensional heavy-tailed distributions.
{\em Ann. Actuar. Scien.}, \textbf{19}, no. 2, 317 -- 349.

\bibitem{konstantinides:passalidis:2023} 
{\sc Konstantinides, D.G., Passalidis, C.D.} \ (2025b)
Background risk model in presence of heavy tails under dependence.
{\em Non. Anal. Mod. Contr.}, \textbf{30}, no. 5, 982--1010.



\bibitem{korshunov:schlegel:schmidt:2003}
{\sc Korshunov, D.A., Schegel, S., Schmidt, F.}\ (2003)
Asymptotic analysis of random walks with dependent heavy-tailed increments.
{\em Siber. Math. J.}, \textbf{44}, no.5, 833--844.


\bibitem{lehmann:1966}
{\sc Lehmann, E.L.}\ (1966)
Some concepts of dependence. 
{\em Ann. Math. Stat.}, \textbf{37}, 1137--1153. 

\bibitem{leipus:siaulys:2020}
{\sc Leipus, R.,  {\v S}iaulys, J.}\ (2020)
On a closure property of convolution equivalent class of distributions.
{\em J. Math. Anal. Appl.}, \textbf{490}, no. 124226.



\bibitem{leipus:siaulys:konstantinides:2023}
{\sc Leipus, R., \v{S}iaulys, J., Konstantinides, D.G.}\ (2023)
{\em Closure Properties for Heavy-Tailed and Related Distributions: An Overview.}
Springer, Cham.



\bibitem{leslie:1989}
{\sc Leslie, J.R.}\ (1989)
On the non-closure under convolution of the subexponential family. 
{\em J. Appl. Probab.}, \textbf{26}, 58--66. 

\bibitem{li:2013}
{\sc Li, J.}\ (2013)
On pairwise quasi-asymptotically independent random variables and  their applications. 
{\em Stat. Prob. Lett.}, \textbf{83}, 2081--2087.

\bibitem{li:2016}
{\sc Li, J.}\ (2016)
Uniform asymptotics for a multi-dimensional time-dependent risk model with multivariate regularly varying claims and stochastic return. 
{\em Insur. Math. Econom.}, \textbf{71}, 195--204.




\bibitem{li:2022a}
{\sc Li, J.}\ (2022)
Asymptotic results on marginal expected shortfalls for dependent risks.
{\em Insur. Math. Econom.}, \textbf{102}, 310--324.


\bibitem{li:sun:2009}
{\sc Li, H., Sun, Y.}\ (2009)
Tail dependence for heavy-tailed scale mixtures of multivariate distributions.
{\em J. Appl. Probab.}, \textbf{46}, 925--937.



\bibitem{li:tang:wu:2010} 
{\sc Li, J., Tang, Q., Wu, R.}\ (2010)
Subexponential tails of discounted aggregate claims in a time-dependent renewal risk model. 
{\em Adv. Appl. Probab.}, \textbf{42}, no. 4, 1126-1146.

\bibitem{liu:shushi:2024}
{\sc Liu, J., Shushi, T.}\ (2024)
Asymptotics of the loss-based tail risk measures in the presence of extreme risks.
{\em Eur. Actuar. J.}, \textbf{14}, 205--224.

\bibitem{liu:woo:2014}
{\sc Liu, J.C., Woo, J.K.}\ (2014)
Asymptotic analysis for risk quantities conditional on ruin for multidimensional heavy-tailed random walks.
{\em Insur. Math. Econom.}, \textbf{55}, 1--9.





\bibitem{meerschaert:scheffler:2001}
{\sc Meerschaert, M.M., Scheffler, H.P.}\ (2001)
{\em Limit distributions for sums of independent random vectors: Heavy tails in theory and practice.}
Wiley, New York.

\bibitem{mikosch:samorodnitsky:2000a} 
{\sc Mikosch, T., Samorodnitsky, G.}\ (2000) 
The supremum of a negative drift random walk with dependent heavy-tailed steps. 
{\em Ann. Appl. Probab.}, \textbf{10}, no. 3, 1025--1064.

\bibitem{mikosch:samorodnitsky:2000b} 
{\sc Mikosch, T., Samorodnitsky, G.}\ (2000) 
Ruin probability with claims modeled by a stationary ergodic stable process. 
{\em Ann. Probab.}, \textbf{28}, no. 4, 1814--1851.

\bibitem{mikosch:wintenberger:2024} 
{\sc Mikosch, T., Wintenberger, O.}\ (2024) 
{\em Extreme Value Theory for Time Series: Models with Power-Law tails}
Springer Nature, Switzerland.


\bibitem{ng:tang:yang:2002}
{\sc Ng, K.W., Tang, Q., Yang, H.}\ (2002)
Maxima of sums of heavy-tailed random variables.
{\em Astin Bull.}, \textbf{32}, no.1, 43--55.


\bibitem{omey:2006}
{\sc Omey, E.} \ (2006)
Subexponential distribution functions in $R^{d}$.
{\em J. Math. Sci.}, \textbf{138}, no. 1, 5434--5449.



\bibitem{pakes:2007}
{\sc Pakes, A.G.} \ (2007)
Convolution equivalence and infinite divisibility: Corrections and corollaries.
{\em J. Appl. Probab.}, \textbf{44}, 295--305.

\bibitem{resnick:2007}
{\sc Resnick, S.}\ (2007) 
{\em Heavy-Tail Phenomena. Probabilistic and Statistical Modeling.} 
Springer, New York.


\bibitem{ross:2014}
{\sc Ross, S.}\ (2014) 
{\em Introduction to Probability Models.} 
Academic Press, Oxford.


\bibitem{samorodnitsky:2016} 
{\sc Samorodnitsky, G.}\ (2016) 
{\em Stochastic Processes and Long Range Dependence.}
Springer, Cham.

\bibitem{samorodnitsky:sun:2016} 
{\sc Samorodnitsky, G., Sun, J.}\ (2016) 
Multivariate subexponential distributions and their applications. 
{\em Extremes}, \textbf{19}, no. 2, 171--196.

\bibitem{samorodnitsky:taqqu:1994} 
{\sc Samorodnitsky, G., Taqqu, M.}\ (1994) 
{\em Stable Non-Gaussian Random Processes: Stochastic Models with Infinite Variance}
Chapman and Hall, New York.








\bibitem{tang:2006}
{\sc Tang, Q.}\ (2006)
The subexponentiality of products revisited.
{\em Extremes}, \textbf{9}, 231--241.



\bibitem{tang:tsitsiashvili:2003a}
{\sc Tang, Q., Tsitsiashvili, G.}\ (2003)
Precise estimates for the ruin probability in finite horizon in a discrete-time model with heavy-tailed insurance and financial risks.
{\em Stoch. Process. Appl.}, \textbf{108}, 299--325.


\bibitem{tang:yuan:2014}
{\sc Tang, Q., Yuan, Z.}\ (2014) 
Randomly weighted sums of subexponential random variables with application to capital allocation.
{\em Extremes}, \textbf{17}, 467--493.



\bibitem{wang:2011}
{\sc Wang, K.}\ (2011)
Randomly weighted sums of dependent subexponential random variables.
{\em Lith. Math. J.}, \textbf{51}, no. 4, 573--586.

\bibitem{wang:su:yang:2024}
{\sc Wang, H., Su, Q., Yang, Y.}\ (2024)
Asymptotics for ruin probabilities in a bidimensional discrete-time risk model with dependence and consistently varying net losses. 
{\em Stochastics}, \textbf{96}, no. 1, 667--695.



\bibitem{watanabe:2008}
{\sc Watanabe, T.}\  (2008) 
Convolution equivalence and distribution of random sums. 
{\em Probab. Theory Relat. Fields}, \textbf{124},  367--397.




\bibitem{xu:foss:wang:2015}
{\sc Xu, H., Foss, S., Wang, Y.}\ (2015)
Convolution and convolution-root properties of long-tailed distributions.
{\em Extremes}, \textbf{18}, 605--628.

\bibitem{xu:cheng:wang:cheng:2018}
{\sc Xu, H., Cheng, F., Wang, Y., Cheng, D.}\ (2018)
A necessary and sufficient condition for the subexponentiality of the product convolution.
{\em Adv. Appl. Probab.}, \textbf{50}, 57--73.














\bibitem{yuan:lu:2023}
{\sc Yuan, M., Lu, D.}\ (2023)
Asymptotics for a time-dependent by-claim model with dependent subexponential claims.
{\em Insur. Math. Econom.}, \textbf{112}, 120--141.

\bibitem{zhu:li:2012}
{\sc Zhu, L., Li, H.}\ (2012) 
Asymptotic analysis of multivariate tail conditional expectations.
{\em N. Amer. Act. J.}, \textbf{16}, 350--363.


\end{thebibliography}
\end{document}